\documentclass[11pt]{article}

\usepackage[a4paper,margin=30mm]{geometry}
\usepackage{amsmath,amssymb,amsthm,mathtools}
\usepackage{bm}
\usepackage{enumitem}
\usepackage{microtype}
\usepackage[
        dvipdfmx,
        hypertexnames=false,
        hyperindex,
        bookmarks=false,
        colorlinks=true,
        linkcolor=magenta,
        citecolor=magenta,
        urlcolor=blue
]{hyperref}
\usepackage{url}
\usepackage{listings}
\lstdefinestyle{magmacode}{
  basicstyle=\footnotesize\ttfamily,
  columns=fullflexible,
  keepspaces=true,
  showstringspaces=false,
  breaklines=true,
    frame=single,
  framesep=4pt
}
\usepackage{color}

\newtheorem{theorem}{Theorem}[section]
\newtheorem{proposition}[theorem]{Proposition}
\newtheorem{lemma}[theorem]{Lemma}
\newtheorem{corollary}[theorem]{Corollary}
\newtheorem{conjecture}[theorem]{Conjecture}

\theoremstyle{definition}
\newtheorem{definition}[theorem]{Definition}

\newtheorem{problem}[theorem]{Problem}
\theoremstyle{remark}
\newtheorem{remark}[theorem]{Remark}

\numberwithin{equation}{section}

\newcommand{\GB}{\operatorname{GB}}
\newcommand{\LM}{\operatorname{LM}}

\newcommand{\ini}{\operatorname{in}}
\newcommand{\gin}{\operatorname{gin}}
\newcommand{\HF}{\operatorname{HF}}
\newcommand{\HS}{\operatorname{HS}}

\newcommand{\Frac}{\operatorname{Frac}}

\title{Specialization of Gr\"obner Bases and\\
a Criterion for Moreno-Soc\'ias' Conjecture\thanks{A preliminary version of some of the results in this paper appeared as part of arXiv:2410.23211.}}

\author{Momonari Kudo\thanks{Fukuoka Institute of Technology, 3-30-1 Wajiro-higashi, Higashi-ku, Fukuoka 811-0295, Japan.  Email: \texttt{m-kudo@fit.ac.jp}.}
\and
Kazuhiro Yokoyama\thanks{Rikkyo University, 3-34-1 Nishi-Ikebukuro, Toshima-ku, Tokyo 171-8501, Japan.  Email: \texttt{kazuhiro-yokoyama@rikkyo.ac.jp}.}}
\date{}

\begin{document}
\maketitle

\vspace{-3mm}

\begin{abstract}
In this paper, we consider homogeneous polynomial sequences with prescribed degrees in a fixed polynomial ring over an infinite field of arbitrary characteristic.
We ask when the leading monomials of the reduced Gr\"obner basis of the ideal generated by the universal generic sequence can be determined from a single specialization of its algebraically independent coefficients.

Our main result is a specialization theorem for reduced Gr\"obner bases up to a prescribed degree over a Noetherian local domain.
If the generic and special fibers have the same Hilbert function through a prescribed degree and the special fiber satisfies a suitable initial-segment condition on its leading monomials, then their reduced Gr\"obner bases agree under specialization through that degree.
Specifically, we work with the graded reverse lexicographic order, for which the weakly reverse lexicographic property of the initial ideal of the special fiber implies the required initial-segment condition.

As a consequence, Moreno-Soc\'ias' conjecture for a fixed number of variables and prescribed degrees admits a one-example criterion:
It is enough to find one specialization whose associated quotient ring has the generic Hilbert series and whose initial ideal is weakly reverse lexicographic.
The Hilbert series predicted by Fr\"oberg's conjecture gives a concrete criterion that verifies Fr\"oberg's and Moreno-Soc\'ias' conjectures simultaneously.
Moreover, an example over a prime field $\mathbb F_p$ satisfying this simultaneous criterion establishes both conjectures for the same type over every field of characteristic zero.
We also combine our criterion with Pardue's reduction to complete intersections, give several computational examples, and establish a corresponding criterion for sequences with a prescribed projective zero.
\end{abstract}

\medskip
\noindent\textit{2020 Mathematics Subject Classification.} Primary 13P10; Secondary 13D40, 13A02.

\smallskip
\noindent\textit{Keywords.} Generic homogeneous ideals, Gr\"obner bases, Hilbert series, Moreno-Soc\'ias' conjecture, Fr\"oberg's conjecture.

\section{Introduction}\label{sec:intro}

Let $K$ be an infinite field and let $R=K[x_1,\ldots,x_n]$ be the polynomial ring in $n$ variables $x_1,\ldots,x_n$ over $K$.
For $d\geq0$, let $R_d$ denote the $K$-vector space of homogeneous polynomials of degree $d$ in $R$.
Fix positive integers $d_1,\ldots,d_m$, put $\tau=(n;d_1,\ldots,d_m)$, $V_\tau:=R_{d_1}\times\cdots\times R_{d_m}$, and, for $\bm F=(f_1,\ldots,f_m)\in V_\tau$, set $I_{\bm F}:=\langle f_1,\ldots,f_m\rangle\subset R$.
We call $V_\tau$ the coefficient space of type $\tau$.
Two classical conjectures predict the behavior of $I_{\bm F}$ for generic $\bm F\in V_\tau$.
Fr\"oberg's conjecture~\cite{Froberg} predicts the Hilbert series of $R/I_{\bm F}$, whereas Moreno-Soc\'ias' conjecture~\cite{MorenoSocias,Pardue} predicts that the initial ideal $\ini_\prec(I_{\bm F})$ is weakly reverse lexicographic when $\prec$ is the graded reverse lexicographic order.
Both conjectures are central problems concerning generic homogeneous ideals and their Gr\"obner bases.

The word ``generic'' is used in two closely related but formally distinct senses.
In Fr\"oberg's original formulation~\cite{Froberg}, the coefficients of the forms are complex numbers that are jointly algebraically independent over $\mathbb Q$.
Here we take the coefficients to be indeterminates over $K$ and consider the resulting forms over the rational function field in these indeterminates.
We call the Hilbert series of the quotient by the ideal generated by these forms the generic Hilbert series of type $\tau$.
In the geometric formulation used, for example, by Pardue~\cite{Pardue}, a property is generic if it holds on a nonempty Zariski-open subset of $V_\tau$.
For properties determined by the initial ideal, these two formulations are equivalent because the reduced Gr\"obner basis is stable under specialization on a suitable nonempty Zariski-open subset.
For Hilbert series, the existence of a nonempty open subset on which the Hilbert series is coefficientwise minimal goes back to Fr\"oberg--L\"ofwall~\cite{FrobergLofwall} and is also discussed by Pardue~\cite{Pardue} and Tani~\cite{Tani}.

For the Hilbert series predicted by Fr\"oberg's conjecture, a stronger one-example principle is known (Proposition~\ref{prop:Frob-one-example}).
For a fixed type $\tau$, the existence of one sequence $\bm F\in V_\tau$ for which $R/I_{\bm F}$ has the predicted Hilbert series forces the quotient associated with the generic sequence to have the same Hilbert series.
This observation is already explicit in Fr\"oberg--Hollman~\cite{FrobergHollman}; see also \cite{Diem}.
By contrast, the existence of a sequence with a weakly reverse lexicographic initial ideal imposes no restriction on the type: such a sequence exists for every type (see Section~\ref{subsec:FMS} below).
The relevant question is therefore whether an appropriate condition on the Hilbert series of $R/I_{\bm F}$ for a single sequence $\bm{F}$, together with the weakly reverse lexicographic property of $\ini_\prec(I_{\bm F})$, is sufficient to force the initial ideal of the generic sequence to be weakly reverse lexicographic.

The main contribution of this paper is a fixed-type criterion formulated in terms of the Hilbert series of the generic quotient, which is more general than the criterion stated in our earlier preprint~\cite{KudoYokoyamaFull}.
We first prove a specialization theorem up to a prescribed degree (Theorem~\ref{thm:fixed-degreewise-specialization}) over a Noetherian local domain.
The local domain need not contain a field, so this formulation includes both specialization of coefficient parameters and reduction modulo a prime.
The theorem requires equality of the Hilbert functions of the generic and special fibers only through that degree, together with an initial-segment condition on the leading monomials of the special fiber in every degree up to the same bound.
The weakly reverse lexicographic property implies this condition.
The proof compares the reduced Gr\"obner bases degree by degree.
Once the leading monomials in lower degrees have been identified, the common Hilbert function determines the number of new basis elements in degree $d$, and the initial-segment condition determines their leading monomials in the special fiber.
We lift the new basis elements of the special fiber and reduce them by the basis elements of the generic fiber already constructed in lower degrees.
In the coefficient matrix of the resulting polynomials, the square submatrix formed by the columns corresponding to these leading monomials reduces to the identity matrix over the residue field.
Its determinant is therefore a unit in the local domain, so row operations over this domain yield the corresponding basis elements of the generic fiber.

Applying the specialization theorem at points of the coefficient space $V_\tau$ yields our main fixed-type criterion (Theorem~\ref{thm:generic-HS-criterion}):
\begin{quote}
For a fixed type $\tau$, Moreno-Soc\'ias' conjecture holds if and only if there exists $\bm F\in V_\tau$ such that $R/I_{\bm F}$ has the generic Hilbert series of type $\tau$ and $\ini_\prec(I_{\bm F})$ is weakly reverse lexicographic.
\end{quote}
Thus the explicit Hilbert series predicted by Fr\"oberg's conjecture is stronger than what is needed for this transfer of initial ideals.
Its advantage is that it gives a directly verifiable condition and establishes Fr\"oberg's and Moreno-Soc\'ias' conjectures simultaneously (Corollary~\ref{cor:Frob-MS}).
The local formulation also gives an arithmetic criterion (Corollary~\ref{cor:arithmetic-lifting}).
If a homogeneous sequence over a finite prime field $\mathbb F_p$ has the Hilbert series predicted by Fr\"oberg's conjecture and a weakly reverse lexicographic initial ideal, then any of its homogeneous lifts over $\mathbb{Z}$ in the same variables and with the same prescribed degrees has the same Hilbert series over $\mathbb Q$, and the reduced Gr\"obner bases over $\mathbb Q$ and $\mathbb F_p$ have the same leading monomials.
Consequently, one such example over $\mathbb F_p$ establishes both conjectures for that type over every field of characteristic zero.
Conversely, if both conjectures hold over $\mathbb Q$ for a fixed type $\tau$, then such an example exists over $\mathbb F_p$ for all but finitely many primes $p$ (Remark~\ref{rem:modular-computation}).
For a fixed homogeneous input with integer coefficients, Corollary~\ref{cor:arithmetic-lifting} also gives a sufficient condition for a prime $p$ to be lucky in modular Gr\"obner basis computation.
This condition can be checked entirely over $\mathbb F_p$.
The assertion about specialization of the reduced Gr\"obner basis can also be deduced from Arnold's results \cite[Theorems~5.6, 5.12, and~6.2]{Arnold}; here it follows directly from our local specialization theorem.

We also record two consequences and extensions.
First, combining our criterion with Pardue's reduction to the case of equally many variables and generators reduces the full Moreno-Soc\'ias conjecture to the construction, for every degree sequence, of one complete intersection with weakly reverse lexicographic initial ideal (Corollary~\ref{cor:CI-construction}).
This connects our criterion with work on generic initial ideals of monomial complete intersections by Cho--Park~\cite{ChoPark} and on generic complete intersections by Capaverde--Gao~\cite{CapaverdeGao}.
Second, the same specialization argument applies to the affine linear parameter space of sequences having a prescribed projective zero (Theorem~\ref{thm:Zp-generic-HS}).
This yields the corresponding one-example criterion in the Krull-dimension-one setting (Theorem~\ref{thm:dimone-criterion}).

A related one-example result for lexicographic order was obtained by Tani~\cite[Theorem~3.5]{Tani}.
Roughly speaking, if the quotient by one ideal has the generic Hilbert series and the ideal has a lexsegment initial ideal, then the initial ideal of the $V$-generic ideals (in Tani's terminology) is lexsegment.
Tani's theorem uses the extremality of lexsegment ideals and exterior powers.
For graded reverse lexicographic order in characteristic zero,
Cho and Park~\cite{ChoPark} formulated a related criterion in terms of the generic initial ideals of two ideals whose quotients have the same Hilbert function.
Our theorem compares the generic and special fibers of a single family in the given coordinates. Under the initial-segment condition on the special fiber, equality of the Hilbert functions through a prescribed degree ensures that the reduced Gr\"obner basis of the generic fiber specializes to that of the special fiber through that degree.
The two fibers may have different characteristics.

The paper is organized as follows.
Section~\ref{sec:prelim} introduces the coefficient spaces, distinguishes the two notions of genericity, and recalls the classical one-example criterion for Fr\"oberg's conjecture.
Section~\ref{sec:specialization} proves the local specialization theorem up to a prescribed degree, derives the fixed-type and arithmetic criteria, gives computational examples, and reduces Moreno-Soc\'ias' conjecture to a problem concerning complete intersections.
Section~\ref{sec:dimone} treats sequences with a prescribed projective zero and the corresponding Krull-dimension-one criterion.
Appendix~\ref{app:marked} provides an alternative, module-theoretic proof using marked families and base change, including the degreewise and prescribed-zero variants.
Appendix~\ref{app:magma} contains \textsc{Magma} code for searching for and verifying the computational examples.

\section{Preliminaries}\label{sec:prelim}

This section fixes notation for homogeneous Gr\"obner bases and weakly reverse lexicographic ideals, and introduces the coefficient spaces and two notions of genericity used throughout the paper.
We also recall Fr\"oberg's and Moreno-Soc\'ias' conjectures, together with the classical one-example criterion for Fr\"oberg's conjecture.

\subsection{Homogeneous Gr\"obner bases and weakly reverse lexicographic ideals}

Let $K$ be a field and let $R=K[x_1,\ldots,x_n]$ with its standard grading, where $\deg x_i=1$ for all $i$.
By a \emph{monomial}, we mean $x^\alpha=x_1^{\alpha_1}\cdots x_n^{\alpha_n}$ for $\alpha=(\alpha_1,\ldots,\alpha_n)\in\mathbb Z_{\geq0}^n$, with coefficient $1$.
Let $\mathbb T$ denote the set of all monomials in $R$.
For a non-negative integer $d \geq 0$, let $R_d$ be the $K$-vector space consisting of $0$ and all homogeneous polynomials of degree $d$, and let $\mathbb T_d$ be the set of monomials of degree $d$.
Then $\mathbb T_d$ is a $K$-linear basis of $R_d$.
We also set $R_e=0$ for $e<0$.
For a subset $S\subseteq R$, let $\langle S\rangle$ denote the ideal of $R$ generated by $S$, with the convention $\langle\emptyset\rangle=0$, the zero ideal.
An ideal generated by a subset of $\mathbb T$ is called a \emph{monomial ideal}.
In particular, the zero ideal is a monomial ideal.

For a finitely generated nonnegatively graded $R$-module $M$, let $M_d$ denote its homogeneous component of degree $d$, so that $M=\bigoplus_{d\geq0}M_d$.
We denote its \emph{Hilbert function} and \emph{Hilbert series} by $\HF_M$ and $\HS_M$, respectively, which are defined by
\[
 \HF_M(d):=\dim_K M_d,
 \qquad
 \HS_M(z):=\sum_{d\geq0}\HF_M(d)z^d.
\]
Unless otherwise specified, $\prec$ denotes the graded reverse lexicographic order with $x_1\succ\cdots\succ x_n$.
For a nonzero polynomial $f\in R$, let $\LM_\prec(f)$ denote its leading monomial with respect to $\prec$.
For an ideal $I\subset R$, its \emph{initial ideal} with respect to $\prec$ is defined by $\ini_\prec(I):=\langle\LM_\prec(f)\mid f\in I\smallsetminus\{0\}\rangle$.
The monomials not belonging to $\ini_\prec(I)$ are called the \emph{standard monomials} of $I$.
For a set $\mathcal G$ of nonzero polynomials in $R$, let $\LM_\prec(\mathcal G):=\{\LM_\prec(g)\mid g\in\mathcal G\}$.
For a homogeneous ideal $I$, let $\GB_\prec(I)$ denote its reduced Gr\"obner basis with respect to $\prec$.
For $d\geq0$, put
\[
 \GB_\prec(I)_d:=\{g\in\GB_\prec(I):\deg g=d\},
 \qquad 
 \GB_\prec(I)_{<d}:=\bigcup_{0\leq e<d}\GB_\prec(I)_e.
\]
We use analogous notation for the set of leading monomials and set $\GB_\prec(0):=\emptyset$.
We use the same notation over any field $E$, with Hilbert functions computed over $E$.
For any coefficient ring $C$, we abbreviate $C[x_1,\ldots,x_n]$ to $C[x]$.

\begin{definition}\label{def:wrl}
A monomial ideal $J\subset R$ is \emph{weakly reverse lexicographic} with respect to $\prec$ if, whenever $t$ is a minimal (monomial) generator of $J$, every monomial $t'$ satisfying $\deg t'=\deg t$ and $t'\succ t$ belongs to $J$.
It is \emph{weakly reverse lexicographic up to degree $D$} if this condition holds for all minimal generators of degree at most $D$.
The zero ideal satisfies these conditions vacuously.
\end{definition}

For a homogeneous ideal $I\subset E[x]$, where $E$ is any field, and $d\geq0$, define
\begin{equation}\label{eq:Bd}
 B_d(I):=
 \left\{
 t\in\mathbb T_d:
 u\nmid t\text{ for every }u\in\LM_\prec(\GB_\prec(I)_{<d})
 \right\}.
\end{equation}
Thus $B_d(I)$ consists of the degree-$d$ monomials that are not divisible by any leading monomial of a reduced Gr\"obner basis element of degree smaller than $d$.
Write $B_d(I)=\{t_1,\ldots,t_w\}$, where $t_1\succ\cdots\succ t_w$.
A subset $S\subseteq B_d(I)$ is called an \emph{initial segment} of $B_d(I)$ if $S=\{t_1,\ldots,t_v\}$ for some integer $0\leq v\leq w$.
Here $v=0$ means $S=\emptyset$.
By the minimality of the leading monomials in a reduced Gr\"obner basis, one has $B_d(I)\cap\ini_\prec(I)=\LM_\prec(\GB_\prec(I)_d)$.
Hence Macaulay's basis theorem~\cite[Theorem~15.3]{Eisenbud} gives
\begin{equation}\label{eq:Bd-HF}
\HF_{E[x]/I}(d)=\#B_d(I)-\#\GB_\prec(I)_d.
\end{equation}

\begin{lemma}\label{lem:wrl-initial-segment}
    Let $E$ be a field, let $I\subset E[x]$ be a homogeneous ideal, and let $D\geq0$ be an integer.
    Then the initial ideal $\ini_\prec(I)$ is weakly reverse lexicographic up to degree $D$ if and only if $\LM_\prec(\GB_\prec(I)_d)$ is an initial segment of $B_d(I)$ for every $0\leq d\leq D$.
\end{lemma}

\begin{proof}
    The forward implication follows from the definition of weakly reverse lexicographic ideals and the equality $B_d(I)\cap\ini_\prec(I)=\LM_\prec(\GB_\prec(I)_d)$.
    Conversely, assume that $\LM_\prec(\GB_\prec(I)_e)$ is an initial segment of $B_e(I)$ for every $0\leq e\leq D$.
    Let $t$ be a minimal monomial generator of $\ini_\prec(I)$ of degree $d\leq D$, and let $t'\succ t$ be another monomial of degree $d$.
    If $t'\notin B_d(I)$, then $t'$ belongs to $\ini_\prec(I)$ by the definition of $B_d(I)$.
    Otherwise, since $\LM_\prec(\GB_\prec(I)_d)$ is an initial segment of $B_d(I)$ containing $t$, we have $t'\in\LM_\prec(\GB_\prec(I)_d)\subset\ini_\prec(I)$.
\end{proof}

We will use this equivalent formulation in the specialization theorem up to a prescribed degree, since it describes the leading monomials of the new reduced Gr\"obner basis elements in each degree.

\begin{remark}
With the polynomial ring $E[x]$ and the monomial order $\prec$ fixed, there is at most one weakly reverse lexicographic monomial ideal $J\subset E[x]$ for which the quotient $E[x]/J$ has a prescribed Hilbert function; see \cite[Section~4, p.~583]{Pardue}.

A corresponding uniqueness statement holds up to any prescribed degree.
Let $D\geq0$ be an integer, and let $I,I'\subset E[x]$ be homogeneous ideals.
Suppose that both $\ini_\prec(I)$ and $\ini_\prec(I')$ are weakly reverse lexicographic up to degree $D$, and that $\HF_{E[x]/I}(d)=\HF_{E[x]/I'}(d)$ for every $0\leq d\leq D$.
Then $\ini_\prec(I)\cap\mathbb T_d=\ini_\prec(I')\cap\mathbb T_d$ for every $0\leq d\leq D$.
Indeed, this follows by induction on $d$.
If the initial ideals agree in lower degrees, then $B_d(I)=B_d(I')$.
By \eqref{eq:Bd-HF}, the sets $\LM_\prec(\GB_\prec(I)_d)$ and $\LM_\prec(\GB_\prec(I')_d)$ have the same cardinality.
By Lemma~\ref{lem:wrl-initial-segment}, both are initial segments of the same set $B_d(I)$, and hence coincide.

Note that this uniqueness statement assumes that {\it both} initial ideals are weakly reverse lexicographic up to degree $D$.
In contrast, Theorem~\ref{thm:fixed-degreewise-specialization} below imposes this condition only on the initial ideal of the special fiber and establishes the corresponding property for the initial ideal of the generic fiber, together with the specialization of reduced Gr\"obner bases.
\end{remark}

\subsection{Coefficient spaces and two notions of genericity}\label{subsec:coeff_sp}

From now on, let $K$ be infinite.
Fix positive integers $d_1,\ldots,d_m$ and put $\tau=(n;d_1,\ldots,d_m)$.
For each $j$, choose a fixed ordering of $\mathbb T_{d_j}$.
Every $f_j\in R_{d_j}$ can be written uniquely as $f_j=\sum_{t\in\mathbb T_{d_j}}a_{j,t}t$ with $a_{j,t}\in K$.
For $\bm F=(f_1,\ldots,f_m)$, we call $\bm a(\bm F):=(a_{j,t})_{1\leq j\leq m,\ t\in\mathbb T_{d_j}}$ its \emph{coefficient vector}, with the entries listed first by $j$ and then in the chosen ordering of $\mathbb T_{d_j}$.
Put
\[
 V_\tau:=R_{d_1}\times\cdots\times R_{d_m}\simeq K^N,
 \quad \text{where} \quad
 N:=\sum_{j=1}^m\binom{n+d_j-1}{d_j}.
\]
Here the isomorphism from $V_{\tau}$ to $K^N$ sends $\bm F$ to its coefficient vector $\bm a(\bm F)$.
We call $\tau$ the \emph{type} of the coefficient space $V_\tau$.
It specifies the ambient number of variables and the prescribed degrees of the entries.
The ambient ring is $R$, even when some variables do not occur in a particular sequence.
We equip $V_\tau$ with the Zariski topology via the identification $V_\tau\simeq K^N$.
This topology is independent of the chosen orderings, since changing them only permutes the coordinates.
For $\bm F=(f_1,\ldots,f_m)\in V_\tau$, set $I_{\bm F}:=\langle f_1,\ldots,f_m\rangle\subset R$.
Let
\[
 \mathcal C_\tau:=\{c_{j,t}:1\leq j\leq m,\ t\in\mathbb T_{d_j}\},
 \qquad
 A_\tau:=K[\mathcal C_\tau],
 \qquad
 L_\tau:=K(\mathcal C_\tau).
\]
Here the $c_{j,t}$ are algebraically independent indeterminates over $K$.
For $1\leq j\leq m$, set
\[
 g_j:=\sum_{t\in\mathbb T_{d_j}}c_{j,t}t\in A_\tau[x_1,\ldots,x_n],
 \qquad
 \bm G_\tau:=(g_1,\ldots,g_m).
\]
We call $\bm G_\tau$ the \emph{universal generic sequence on $V_\tau$}.
Here ``universal'' refers to the following property: For every commutative $K$-algebra $C$, each sequence in $\prod_{j=1}^m C[x]_{d_j}$ is obtained from $\bm G_\tau$ by a unique $K$-algebra homomorphism $A_\tau\to C$ applied coefficientwise.
We write
\[
 I_\tau^{\mathrm{gen}}:=\langle\bm G_\tau\rangle\subset L_\tau[x_1,\ldots,x_n].
\]
For a fixed monomial order $\prec$, we call $\ini_\prec(I_\tau^{\mathrm{gen}})$ the \emph{initial ideal of the universal generic sequence} on $V_\tau$.
It should not be confused with the generic initial ideal $\gin_\prec(I)$ obtained from a general linear change of coordinates of a fixed ideal $I$.
For $\bm F\in V_\tau$ with coefficient vector $(a_{j,t})_{j,t}$, let $\pi_{\bm F}:A_\tau\to K$ be the coefficient specialization $c_{j,t}\mapsto a_{j,t}$, and put
\[
 \mathfrak m_{\bm F}:=\ker\pi_{\bm F}
 =\langle c_{j,t}-a_{j,t}\mid 1\leq j\leq m,\ t\in\mathbb T_{d_j} \rangle,
 \qquad
 \mathcal O_{\bm F}:=(A_\tau)_{\mathfrak m_{\bm F}}\subset L_\tau.
\]
The map $\pi_{\bm F}$ extends naturally to the residue map $\mathcal O_{\bm F}\to K$, given by $a/b\mapsto\pi_{\bm F}(a)/\pi_{\bm F}(b)$ for $a,b\in A_\tau$ with $b\notin\mathfrak m_{\bm F}$ (i.e., $\pi_{\bm F}(b)\ne0$).
Its kernel is $\mathfrak m_{\bm F}\mathcal O_{\bm F}$, the unique maximal ideal of $\mathcal O_{\bm F}$.
We also write $\pi_{\bm F}$ for its coefficientwise extension $\mathcal O_{\bm F}[x]\to K[x]=R$; in particular, $\pi_{\bm F}(\bm G_\tau)=\bm F$.
A rational function in $L_\tau$ is said to be \emph{regular at $\bm F$} if it belongs to $\mathcal O_{\bm F}$.
Thus a rational function can be specialized at $\bm F$ if and only if it can be written as $a/b$ with $a,b\in A_\tau$ and $b\notin\mathfrak m_{\bm F}$ (i.e., $\pi_{\bm F}(b)\ne0$).

\begin{definition}\label{def:generic-property}
A property $\mathcal P$ of sequences in $V_\tau$ is \emph{generic} if there exists a nonempty Zariski-open subset $U\subset V_\tau$ such that every $\bm F\in U$ satisfies $\mathcal P$.
The same definition applies to any $K$-linear subspace $Z\subset V_\tau$, using the induced Zariski topology on $Z$.
\end{definition}

\begin{remark}[The two formulations]\label{rem:two-genericities}
The assertion that $\bm G_\tau$ satisfies a property is understood over the rational function field $L_\tau$.
In contrast, Definition~\ref{def:generic-property} concerns a property of sequences over $K$ that holds on a nonempty Zariski-open subset of $V_\tau$.
These are distinct formulations.
For properties determined solely by the monomials contained in the initial ideal, Proposition~\ref{prop:initial-property} below shows that they are equivalent.

In his original formulation over $\mathbb C$, Fr\"oberg required the coefficients of all forms to be jointly algebraically independent over $\mathbb Q$ \cite[pp.~118--119]{Froberg}.
Here we allow an arbitrary infinite field $K$ and introduce indeterminates that are algebraically independent over $K$.
Let $k_0$ be the prime field of $K$.
Since $g_j\in k_0(\mathcal C_\tau)[x]$, the reduced Gr\"obner basis of $\langle g_1,\ldots,g_m\rangle$ over $k_0(\mathcal C_\tau)$ remains the reduced Gr\"obner basis after extending the coefficient field to $L_\tau=K(\mathcal C_\tau)$.
Hence the leading monomials and the Hilbert series of the quotient are unchanged by this field extension.
The corresponding invariance of the Hilbert series is also stated in \cite[Lemma and Definition~4]{Diem}.

In characteristic zero, let $\bm a=(a_{j,t})_{j,t}$ be a tuple of complex numbers that are jointly algebraically independent over $\mathbb Q$, and set $f_j=\sum_{t\in\mathbb T_{d_j}}a_{j,t}t$.
The assignment $c_{j,t}\mapsto a_{j,t}$ induces a field isomorphism $\mathbb Q(\mathcal C_\tau)\simeq\mathbb Q(\bm a)$.
Since reduced Gr\"obner bases are preserved under field isomorphisms and extension of the coefficient field, the reduced Gr\"obner bases of $\langle f_1,\ldots,f_m\rangle\subset\mathbb C[x]$ and $I_\tau^{\mathrm{gen}}\subset L_\tau[x]$ have the same leading monomials, and their quotient rings have the same Hilbert series.
\end{remark}

\begin{proposition}[Generic stability of a reduced Gr\"obner basis]\label{prop:generic-stability}
Fix a monomial order $\prec$, and let $\mathcal G_\tau$ be the reduced Gr\"obner basis of $I_\tau^{\mathrm{gen}}$ in $L_\tau[x]$ with respect to $\prec$.
There exists a nonempty Zariski-open subset $U\subset V_\tau$ such that, for every $\bm F\in U$,
\begin{enumerate}[label=\textup{(\roman*)}]
 \item all coefficients of the elements of $\mathcal G_\tau$ belong to $\mathcal O_{\bm F}$;
 \item $\pi_{\bm F}(\mathcal G_\tau)=\GB_\prec(I_{\bm F})$;
 \item $\LM_\prec(\mathcal G_\tau)=\LM_\prec(\GB_\prec(I_{\bm F}))$.
\end{enumerate}
\end{proposition}

\begin{proof}
This is the standard stability of a Gr\"obner basis under specialization; see, for example, \cite{Kalkbrener,Nabeshima}.
Write $\mathcal G_\tau=\{h_1,\ldots,h_s\}$.
Since $\langle h_1,\ldots,h_s\rangle=\langle g_1,\ldots,g_m\rangle$ in $L_\tau[x]$, there exist $u_{ij},v_{ji}\in L_\tau[x]$ such that $h_i=\sum_{j=1}^m u_{ij}g_j$ and $g_j=\sum_{i=1}^s v_{ji}h_i$.
Let $\Delta_\tau\in A_\tau\smallsetminus\{0\}$ be the product of the nonzero denominators occurring in the coefficients of the $h_i$, $u_{ij}$, and $v_{ji}$.
All these polynomials then belong to $B[x]$ for $B:=A_\tau[\Delta_\tau^{-1}]$.
The identities above give $\langle h_1,\ldots,h_s\rangle_{B[x]}=\langle g_1,\ldots,g_m\rangle_{B[x]}$.
The leading monomial of every nonzero polynomial in this ideal is divisible by some $\LM_\prec(h_i)$, since $\mathcal G_\tau$ is a Gr\"obner basis over $L_\tau$.
As each $h_i$ is monic, $\mathcal G_\tau$ is also a Gr\"obner basis of this ideal over $B$.

Let $\bm F\in V_\tau$ satisfy $\pi_{\bm F}(\Delta_\tau)\ne0$.
Then $B\subseteq\mathcal O_{\bm F}$, so that $\mathcal G_\tau\subset\mathcal O_{\bm F}[x]$.
Since each $h_i$ is monic, its leading monomial is unchanged under specialization.
Applying \cite[Theorem~3.1]{Kalkbrener}, with all leading coefficients equal to $1$, shows that $\pi_{\bm F}(\mathcal G_\tau)$ is a Gr\"obner basis of $I_{\bm F}$.
Its reducedness is also preserved, since its leading monomials are unchanged and no new monomials occur.
Hence $\pi_{\bm F}(\mathcal G_\tau)=\GB_\prec(I_{\bm F})$ and $\LM_\prec(\mathcal G_\tau)=\LM_\prec(\GB_\prec(I_{\bm F}))$.
Finally, since $K$ is infinite and $\Delta_\tau\ne0$, the set $U:=\{\bm F\in V_\tau\mid\pi_{\bm F}(\Delta_\tau)\ne0\}$ is a nonempty Zariski-open subset of $V_\tau$.
\end{proof}

The same proof applies to any finite polynomial family with coefficients in a polynomial ring over $K$; it does not require that all coefficients of the input polynomials be independent parameters.

\begin{definition}\label{def:generic-HS}
The \emph{generic Hilbert series of type $\tau$} is
\[
 H_\tau^{\mathrm{gen}}(z):=
 \HS_{L_\tau[x_1,\ldots,x_n]/I_\tau^{\mathrm{gen}}}(z).
\]
We also call $\HF_{L_\tau[x]/I_\tau^{\mathrm{gen}}}$ the \emph{generic Hilbert function of type $\tau$}.
\end{definition}

\begin{corollary}\label{cor:generic-HS-open}
The equality $\HS_{R/I_{\bm F}}(z)=H_\tau^{\mathrm{gen}}(z)$ holds for $\bm{F}$ in a nonempty Zariski-open subset $U$ of $V_\tau$.
Moreover, this open subset can be chosen so that the initial ideals of $I_{\bm F}$ and $I_\tau^{\mathrm{gen}}$ have the same minimal monomial generators for every $\bm F$ in $U$.
\end{corollary}

\begin{proof}
Apply Proposition~\ref{prop:generic-stability} and Macaulay's basis theorem~\cite[Theorem~15.3]{Eisenbud}.
\end{proof}

For every $\bm F\in V_\tau$ and $d\geq0$, the generic Hilbert function satisfies the inequality $\HF_{L_\tau[x]/I_\tau^{\mathrm{gen}}}(d)\leq\HF_{R/I_{\bm F}}(d)$.
Indeed, in each degree, the rank of the coefficient matrix of the monomial multiples of $g_1,\ldots,g_m$ cannot increase under specialization; see also \cite[Appendix~B, Lemma~3]{Diem}.
The same rank argument applies to a homogeneous family whose coefficients lie in a polynomial ring over $K$.
When all $d_j$ are at least $2$, this coefficientwise minimality is also stated explicitly in \cite[Proposition~2.3]{FrobergLofwall2026}.

\begin{proposition}\label{prop:initial-property}
Let $\mathcal Q$ be a property determined solely by the set of monomials contained in a monomial ideal.
The following are equivalent:
\begin{enumerate}[label=\textup{(\roman*)}]
 \item $\ini_\prec(I_\tau^{\mathrm{gen}})$ satisfies $\mathcal Q$;
 \item the property that $\ini_\prec(I_{\bm F})$ satisfies $\mathcal Q$ is generic on $V_\tau$.
\end{enumerate}
\end{proposition}

\begin{proof}
The implication \textup{(i)} $\Rightarrow$ \textup{(ii)} follows from Proposition~\ref{prop:generic-stability}.
Conversely, let $U$ be a nonempty open subset as in \textup{(ii)}, and let $U'$ be the nonempty open subset supplied by Proposition~\ref{prop:generic-stability}.
Since $V_\tau$ is irreducible, we may choose $\bm F\in U\cap U'$.
Then $\ini_\prec(I_{\bm F})$ satisfies $\mathcal Q$ and contains the same monomials as $\ini_\prec(I_\tau^{\mathrm{gen}})$.
Thus $\ini_\prec(I_\tau^{\mathrm{gen}})$ satisfies $\mathcal Q$.
\end{proof}

\subsection{Fr\"oberg's and Moreno-Soc\'ias' conjectures}\label{subsec:FMS}

For a formal power series $h(z)=\sum_{d\geq0}a_dz^d$ with integer coefficients in the variable $z$, let $[h(z)]$ denote the series obtained by truncating $h(z)$ immediately before its first nonpositive coefficient.
The coefficient of $[h(z)]$ in degree $d$ is $a_d$ if $a_0,\ldots,a_d$ are all positive, and is zero otherwise.
For two power series $A(z)$ and $B(z)$, write $A\succeq_{\mathrm{lex}}B$ if $A=B$ or, at the first degree in which their coefficients differ, the coefficient of $A$ is larger.
Put
\begin{equation}\label{eq:Frob-series}
 H_\tau^{\mathrm{Fr}}(z):=
 \left[
 \frac{\prod_{j=1}^m(1-z^{d_j})}{(1-z)^n}
 \right].
\end{equation}
We use $H_\tau^{\mathrm{Fr}}(z)$ to denote the Hilbert series predicted by Fr\"oberg's conjecture for type $\tau$.

\begin{conjecture}[Fr\"oberg]\label{conj:Froberg}
For every type $\tau$, one has $H_\tau^{\mathrm{gen}}(z)=H_\tau^{\mathrm{Fr}}(z)$.
\end{conjecture}

Fr\"oberg's argument~\cite[pp.~121--123]{Froberg} gives $H_\tau^{\mathrm{gen}}(z)\succeq_{\mathrm{lex}}H_\tau^{\mathrm{Fr}}(z)$:
Although the original proof is given over $\mathbb C$ for forms of degrees greater than one, the argument using multiplication by the successive forms applies over any field, including finite fields, and for arbitrary positive degrees.
By the coefficientwise minimality of the generic Hilbert series shown in Section~\ref{subsec:coeff_sp}, it follows that $\HS_{R/I_{\bm F}}(z)\succeq_{\mathrm{lex}}H_\tau^{\mathrm{Fr}}(z)$ for every $\bm F\in V_\tau$.

The following fixed-type criterion is classical.
The implication from one example to the generic sequence follows from Fr\"oberg's lower bound and minimality of the generic Hilbert series, and was explicitly used by Fr\"oberg--Hollman \cite[p.~150]{FrobergHollman}; see also \cite{Froberg,Diem}.

\begin{proposition}[Classical one-example criterion for Fr\"{o}berg's conjecture]\label{prop:Frob-one-example}
For a fixed type $\tau=(n;d_1,\ldots,d_m)$, the following are equivalent:
\begin{enumerate}[label=\textup{(\roman*)}]
 \item $H_\tau^{\mathrm{gen}}(z)=H_\tau^{\mathrm{Fr}}(z)$;
 \item the equality $\HS_{R/I_{\bm F}}(z)=H_\tau^{\mathrm{Fr}}(z)$ holds for any $\bm{F}$ in a nonempty Zariski-open subset of $V_\tau$;
 \item there exists one $\bm F\in V_\tau$ such that $\HS_{R/I_{\bm F}}(z)=H_\tau^{\mathrm{Fr}}(z)$.
\end{enumerate}
\end{proposition}

\begin{proof}
The implication (i) $\Rightarrow$ (ii) follows from Corollary~\ref{cor:generic-HS-open}, while (ii) $\Rightarrow$ (iii) is immediate.
Assume (iii).
By the coefficientwise minimality of the generic Hilbert series established above, the coefficients of $H_\tau^{\mathrm{gen}}(z)$ are not larger than those of $\HS_{R/I_{\bm F}}(z)=H_\tau^{\mathrm{Fr}}(z)$.
On the other hand, Fr\"oberg's inequality gives $H_\tau^{\mathrm{gen}}(z)\succeq_{\mathrm{lex}}H_\tau^{\mathrm{Fr}}(z)$.
If these two series differed, their first differing coefficient would have to be both smaller and larger in $H_\tau^{\mathrm{gen}}(z)$.
Thus they are equal, proving (i).
\end{proof}

\begin{conjecture}[Moreno-Soc\'ias]\label{conj:MS}
For every type $\tau$, the initial ideal $\ini_\prec(I_\tau^{\mathrm{gen}})$ of the ideal $I_\tau^{\mathrm{gen}}$ generated by the universal generic sequence $\bm{G}_{\tau}$ is weakly reverse lexicographic.
\end{conjecture}

By Proposition~\ref{prop:initial-property}, Conjecture~\ref{conj:MS} is equivalent to the open-set formulation for each fixed type.
In contrast, the \emph{existence} of a specialization with weakly reverse lexicographic initial ideal alone imposes no restriction on the type:
For example, $\bm F=(x_1^{d_1},\ldots,x_1^{d_m})$ generates the weakly reverse lexicographic monomial ideal $\langle x_1^{\min_j d_j}\rangle$ for every type.
Pardue proved that Moreno-Soc\'ias' conjecture for all types implies Fr\"oberg's conjecture for all types \cite[Theorem~2]{Pardue}.
The criteria in Section~\ref{sec:specialization}, by contrast, are stated separately for each fixed type.

\section{Specialization up to a prescribed degree and applications}\label{sec:specialization}

In this section, we first establish a specialization theorem for reduced Gr\"obner bases up to a prescribed degree over a Noetherian local domain.
The theorem requires equality of the Hilbert functions of the generic and special fibers through that degree, but not an explicit formula for these functions.
We then derive fixed-type criteria for Fr\"oberg's and Moreno-Soc\'ias' conjectures, as well as a criterion for lifting from a finite prime field to characteristic zero.
The local domain need not contain a field, and its residue field may be finite, so this formulation applies both to specialization of coefficient parameters and to reduction modulo a prime.
For applications to coefficient spaces, we retain the notation of Section~\ref{sec:prelim}, in particular the infinite field $K$, the polynomial ring $R=K[x]=K[x_1,\ldots,x_n]$, the type $\tau=(n;d_1,\ldots,d_m)$, and the coefficient space $V_\tau=R_{d_1}\times \cdots \times R_{d_m}$.

\subsection{Specialization up to a prescribed degree}

The following theorem compares the reduced Gr\"obner bases of the generic and special fibers by combining equality of Hilbert functions with the initial-segment property in Lemma~\ref{lem:wrl-initial-segment}.

\begin{theorem}[Specialization over a local domain]\label{thm:fixed-degreewise-specialization}
Let $(B,\mathfrak m)$ be a Noetherian local domain, and put $L:=\operatorname{Frac}(B)$ and $k:=B/\mathfrak m$.
Let $\pi_k:B[x]\to k[x]$ be the coefficientwise residue homomorphism.
For a homogeneous ideal $I_B\subset B[x]$, set $I_L:=I_BL[x]$ and $I_k:=\pi_k(I_B)\subset k[x]$.
Since $\pi_k$ is a surjective graded homomorphism, $I_k$ is a homogeneous ideal of $k[x]$.
Let $D$ be a nonnegative integer, and assume that
\begin{equation}\label{eq:fixed-HF-through-D}
\HF_{k[x]/I_k}(d)=\HF_{L[x]/I_L}(d)\qquad(0\leq d\leq D).
\end{equation}
Assume moreover that $\ini_\prec(I_k)$ is weakly reverse lexicographic up to degree $D$.
Equivalently, by Lemma~\ref{lem:wrl-initial-segment}, the set $\LM_\prec\bigl(\GB_\prec(I_k)_d\bigr)$ is an initial segment of $B_d(I_k)$ for every $0\leq d\leq D$.
Then, for every $0\leq d\leq D$, one has $\GB_\prec(I_L)_d\subset B[x]$ and
\begin{equation}\label{eq:fixed-GB-specialization}
\pi_k\bigl(\GB_\prec(I_L)_d\bigr)=\GB_\prec(I_k)_d.
\end{equation}
In particular,
\begin{equation}\label{eq:fixed-LM-specialization}
\LM_\prec\bigl(\GB_\prec(I_L)_d\bigr)=\LM_\prec\bigl(\GB_\prec(I_k)_d\bigr) \qquad(0\leq d\leq D).
\end{equation}
Here monomials over $L$ and $k$ are identified by their exponent vectors.
\end{theorem}

\begin{proof}
Put $\mathcal G:=\GB_\prec(I_L)$ and $\mathcal H:=\GB_\prec(I_k)$.
Since every element of a reduced Gr\"obner basis is monic, specialization preserves its leading monomial whenever all its coefficients belong to $B$.
Thus \eqref{eq:fixed-LM-specialization} follows from \eqref{eq:fixed-GB-specialization}.

We prove $\mathcal G_d\subset B[x]$ and \eqref{eq:fixed-GB-specialization} by induction on $d$.
In degree $0$, equality of Hilbert functions implies that either both degree-zero parts of the reduced Gr\"obner bases are empty or both are $\{1\}$, so the assertion holds.
Suppose that it has been proved in every degree smaller than
$d\leq D$.
Then $\LM_\prec(\mathcal G_{<d}) = \LM_\prec(\mathcal H_{<d})$, and hence $B_d(I_L) = B_d(I_k)$.
Write this common set as $B_d:=B_d(I_L) = B_d(I_k)=\{t_1\succ\cdots\succ t_w\}$, and put $v:=\#\mathcal H_d$ and $\bar v:=\#\mathcal G_d$.

\medskip
\noindent \textbf{Step 1: The number and leading monomials of new basis elements.}
By~\eqref{eq:Bd-HF}, we have $\HF_{k[x]/I_k}(d)=w-v$ and $\HF_{L[x]/I_L}(d)=w-\bar v$.
The equality \eqref{eq:fixed-HF-through-D} therefore gives $v=\bar v$.
If $v=0$, there is nothing to prove in degree $d$.
Hence assume $v>0$.
By Lemma~\ref{lem:wrl-initial-segment}, the weakly reverse lexicographic hypothesis implies that $\LM_\prec(\mathcal H_d)=\{t_1,\ldots,t_v\}$, and we may order $\mathcal H_d=\{h_1,\ldots,h_v\}$ so that
\begin{equation}\label{eq:fixed-top-v}
\LM_\prec(h_i)=t_i \qquad(1\leq i\leq v).
\end{equation}
We now construct polynomials in $I_L$ that specialize to $h_1,\ldots,h_v$ and have leading monomials $t_1,\ldots,t_v$.
The equality $v=\#\mathcal G_d$ will then allow us to identify these polynomials with $\mathcal G_d$.

\medskip
\noindent \textbf{Step 2: Lifting and reduction by lower-degree basis elements.}
For each $1\leq i\leq v$, choose a homogeneous polynomial $\ell_i^{(0)}\in(I_B)_d$ with $I_B\subset I_L\cap B[x]$ such that $\pi_k(\ell_i^{(0)})=h_i$.
Such a lift exists because $I_k=\pi_k(I_B)$ and $I_B$ is homogeneous.

Using the elements of $\mathcal G_{<d}$, we successively reduce the terms of $\ell_i^{(0)}$ whose monomials do not belong to $B_d$.
By the induction hypothesis, every element of $\mathcal G_{<d}$ has coefficients in $B$ and specializes to the corresponding element of $\mathcal H_{<d}$, with the same leading monomial.
At every reduction step, suppose that the current polynomial is $p\in B[x]$ with $\pi_k(p)=h_i$, and let $t\in\operatorname{Supp}(p)\smallsetminus B_d$ be the largest monomial to be reduced.
By the definition of $B_d$, there exists $g\in\mathcal G_{<d}$ such that $\LM_\prec(g)\mid t$.
Since $h_i$ is reduced with respect to $\mathcal H_{<d}$, the coefficient of $t$ in $h_i$ is zero.
Hence the coefficient of $t$ in $p$ belongs to $\mathfrak m$ and specializes to zero.
It follows that the reduction does not change the specialization.
Moreover, all new monomials introduced by the reduction are smaller than $t$.
Since the largest monomial to be reduced strictly decreases at each step and $\prec$ is a well-order, the procedure terminates.
Denote the resulting polynomial by $\ell_i$.
Then we have $\ell_i\in I_L \cap B[x]$ and
\begin{equation}\label{eq:fixed-lifts}
\pi_k(\ell_i)=h_i, \qquad \operatorname{Supp}(\ell_i)\subseteq B_d.
\end{equation}

\medskip
\noindent \textbf{Step 3: Row reduction and identification
with the generic basis.}
We next use row operations to obtain monic polynomials with leading monomials $t_1,\ldots,t_v$, without changing their specializations.
Let $M\in\operatorname{Mat}_{v\times w}(B)$ be the coefficient matrix of $\ell_1,\ldots,\ell_v$ with respect to the ordered monomials $t_1,\ldots,t_w$, and let $M_0\in\operatorname{Mat}_v(B)$ be its submatrix consisting of the first $v$ columns.
Here $\pi_k$ is applied entrywise to matrices over $B$.
Since the $h_i$ belong to the reduced Gr\"obner basis $\mathcal{H}$, \eqref{eq:fixed-top-v} and \eqref{eq:fixed-lifts} give $\pi_k(M_0)=I_v$, whence $\pi_k(\det(M_0)) = \det(\pi_k(M_0)) = 1$.
Hence $\det(M_0)$ does not belong to the maximal ideal $\mathfrak{m}$, so $\det(M_0)$ is a unit in $B$.
Therefore, $M_0$ is invertible over $B$, i.e., $M_0^{-1}\in \mathrm{Mat}_{v\times v}(B)$.

Let $\overline h_1,\ldots,\overline h_v$ be the polynomials corresponding to the rows of the matrix $M_0^{-1}M$ over $B$, so that $\overline h_i \in I_L \cap B[x_1,\ldots,x_n]$.
Since $\pi_k(M_0^{-1}M) = \pi_k(M_0)^{-1}\pi_k(M) = \pi_k(M)$, we have $\pi_k(\overline h_i)=h_i$ for $1\leq i \leq v$.
Here, $M_0^{-1}M$ has the form $(I_v\mid *)$, so each $\overline h_i$ is monic with $\LM_\prec(\overline h_i)=t_i$.
As $\overline h_i\in I_L$ and $t_i\in B_d$, the monomials $t_1,\ldots,t_v$ are distinct minimal generators of $\ini_\prec(I_L)$.
Since $\#\mathcal G_d=v$, we obtain $\LM_\prec(\mathcal G_d)=\{t_1,\ldots,t_v\}$.
Moreover, every nonleading monomial occurring in an $\overline h_i$ belongs to $\{t_{v+1},\ldots,t_w\}$, which consists of
standard monomials of $I_L$.
Therefore, by the uniqueness of the reduced Gr\"obner basis, we have $\mathcal G_d=\{\overline h_1,\ldots,\overline h_v\}$, which proves \eqref{eq:fixed-GB-specialization} in degree $d$ and completes the induction.
\end{proof}

\begin{corollary}[Specialization of the full reduced Gr\"obner basis]\label{cor:fixed-full-specialization}
Let $\bm F\in V_\tau$ satisfy $\HS_{R/I_{\bm F}}(z) = H_\tau^{\mathrm{gen}}(z)$.
If $\ini_\prec(I_{\bm F})$ is weakly reverse lexicographic, then $\GB_\prec(I_\tau^{\mathrm{gen}})\subset\mathcal O_{\bm F}[x]$ and
\[
\pi_{\bm F} \bigl( \GB_\prec(I_\tau^{\mathrm{gen}}) \bigr) = \GB_\prec(I_{\bm F}),
\qquad
\LM_\prec \bigl( \GB_\prec(I_\tau^{\mathrm{gen}}) \bigr)
= \LM_\prec \bigl( \GB_\prec(I_{\bm F}) \bigr).
\]
In particular, the two initial ideals have the same minimal monomial generators, and hence $\ini_\prec(I_\tau^{\mathrm{gen}})$ is weakly reverse lexicographic.
\end{corollary}

\begin{proof}
Apply Theorem~\ref{thm:fixed-degreewise-specialization} with $B:=\mathcal O_{\bm F}$ and $I_B:=\langle g_1,\ldots,g_m\rangle\subset B[x]$.
Its fraction field and residue field are $L_\tau$ and $K$, and the ideals $I_L$ and $I_k$ in that theorem are $I_\tau^{\mathrm{gen}}$ and $I_{\bm F}$, respectively.
Choose $D\geq0$ at least as large as the degree of every element of the two reduced Gr\"obner bases.
Since $\ini_\prec(I_{\bm F})$ is weakly reverse lexicographic, the second hypothesis of Theorem~\ref{thm:fixed-degreewise-specialization} holds in every degree.
The Hilbert-function hypothesis follows from the assumed equality of Hilbert series.
Applying the theorem gives the assertion.
\end{proof}

For a prime $p$, put $\mathbb Z_{(p)}:=\{a/b\mid a,b\in\mathbb Z,\ p\nmid b\}$.
Theorem~\ref{thm:fixed-degreewise-specialization} also applies to reduction modulo $p$ by taking $B=\mathbb Z_{(p)}$, whose fraction field is $\mathbb Q$ and whose residue field is $\mathbb F_p$.
If the Hilbert functions of the generic and special fibers agree in every degree and the initial ideal of the special fiber is weakly reverse lexicographic, the theorem gives specialization of the full reduced Gr\"obner basis.
In Section~\ref{sec:fixedtype}, we give a condition on the special fiber that guarantees the required equality of Hilbert functions.

\subsection{Fixed-type criteria and computational examples}\label{sec:fixedtype}

We first return to the coefficient space $V_\tau$ over $K$ and derive the fixed-type criteria.
We then turn to reduction modulo a prime and computational examples.

For $\bm F\in V_\tau$, we define the properties $\mathcal H_{\tau}^{\mathrm{gen}}(\bm{F})$, $\mathcal H_{\tau}^{\mathrm{Fr}}(\bm{F})$, and $\mathcal W({\bm F})$ as follows:
\[
\begin{aligned}
\mathcal H_\tau^{\mathrm{gen}}(\bm F): & \ \HS_{R/I_{\bm F}}(z)=H_\tau^{\mathrm{gen}}(z),\qquad
\mathcal H_\tau^{\mathrm{Fr}}(\bm F) : \ \HS_{R/I_{\bm F}}(z)=H_\tau^{\mathrm{Fr}}(z),\\
\mathcal W({\bm F}):& \ \ini_\prec(I_{\bm F}) \text{ is weakly reverse lexicographic}.
\end{aligned}
\]
By Corollary~\ref{cor:generic-HS-open}, the property $\mathcal H_{\tau}^{\mathrm{gen}}(\bm{F})$ holds for every $\bm{F}$ in a nonempty Zariski-open subset of $V_\tau$.

\begin{theorem}[One-example criterion for Moreno-Soc\'ias' conjecture]\label{thm:generic-HS-criterion}
For every fixed type $\tau$, the following are equivalent:
\begin{enumerate}[label=\textup{(\roman*)}]
 \item $\ini_\prec(I_\tau^{\mathrm{gen}})$ is weakly reverse
 lexicographic;
 \item $\mathcal H_{\tau}^{\mathrm{gen}}({\bm F})\wedge\mathcal W({\bm F})$ holds for every $\bm{F}$ in a nonempty Zariski-open subset of $V_\tau$;
 \item there exists one $\bm F\in V_\tau$ satisfying $\mathcal H_{\tau}^{\mathrm{gen}}({\bm F})\wedge\mathcal W({\bm F})$.
\end{enumerate}
\end{theorem}

\begin{proof}
The implication (i) $\Rightarrow$ (ii) follows immediately from Proposition~\ref{prop:generic-stability}, and the implication (ii) $\Rightarrow$ (iii) is trivial.
The implication (iii) $\Rightarrow$ (i) follows from Corollary~\ref{cor:fixed-full-specialization}.
\end{proof}

\begin{remark}
Since $\mathcal H_{\tau}^{\mathrm{gen}}({\bm F})$ holds for every $\bm{F}$ in a nonempty Zariski-open subset of $V_\tau$, condition (ii) of Theorem~\ref{thm:generic-HS-criterion} is equivalent to the usual generic formulation of Moreno-Soc\'ias' conjecture for the fixed type $\tau$.
Thus the theorem may be read as
\[
 \text{Moreno-Soc\'ias for }\tau
 \quad\Longleftrightarrow\quad
 \exists\,\bm F\in V_\tau \text{ such that }
 \mathcal H_{\tau}^{\mathrm{gen}}({\bm F})\wedge\mathcal W({\bm F}) \text{ holds}.
\]
\end{remark}

Combining Theorem~\ref{thm:generic-HS-criterion} with Proposition~\ref{prop:Frob-one-example} gives the following simultaneous criterion:

\begin{corollary}[Simultaneous criterion for the two conjectures]\label{cor:Frob-MS}
For every fixed type $\tau$, the following are equivalent:
\begin{enumerate}[label=\textup{(\roman*)}]
 \item $H_\tau^{\mathrm{gen}}(z)=H_\tau^{\mathrm{Fr}}(z)$ and $\ini_\prec(I_\tau^{\mathrm{gen}})$ is weakly reverse lexicographic;
 \item $\mathcal H_{\tau}^{\mathrm{Fr}}({\bm F})\wedge\mathcal W({\bm F})$ holds for every $\bm{F}$ in a nonempty Zariski-open subset of $V_\tau$;
 \item there exists one $\bm F\in V_\tau$ satisfying $\mathcal H_{\tau}^{\mathrm{Fr}}({\bm F})\wedge\mathcal W({\bm F})$.
\end{enumerate}
\end{corollary}

\begin{proof}
The implication (i) $\Rightarrow$ (ii) follows from Proposition \ref{prop:generic-stability}, and (ii) $\Rightarrow$ (iii) is immediate.
Assume (iii).
By Proposition~\ref{prop:Frob-one-example}, we have $H_\tau^{\mathrm{gen}}(z)=H_\tau^{\mathrm{Fr}}(z)$.
Thus the chosen $\bm F$ satisfies $\mathcal H_{\tau}^{\mathrm{gen}}({\bm F})\wedge\mathcal W({\bm F})$, and Theorem~\ref{thm:generic-HS-criterion} proves (i).
\end{proof}

The local formulation also gives a criterion using an example over a finite prime field.
The assertions about Hilbert series and reduced Gr\"obner bases in Corollary~\ref{cor:arithmetic-lifting} can also be deduced from Lemma~\ref{lem:hilbert-series-lifting} and Arnold's results \cite{Arnold}; see Remark~\ref{rem:modular-computation}.
Here we obtain the assertions about reduced Gr\"obner bases from Theorem~\ref{thm:fixed-degreewise-specialization}.
We begin with the following comparison of Hilbert functions.

\begin{lemma}[Hilbert functions under reduction modulo a prime]\label{lem:hilbert-series-lifting}
Let $p$ be a prime and let $\tau=(n;d_1,\ldots,d_m)$.
For $1\leq j\leq m$, let $f_j\in\mathbb Z[x]_{d_j}$ and write $\overline f_j\in\mathbb F_p[x]_{d_j}$ for its reduction modulo $p$.
Put $I:=\langle f_1,\ldots,f_m\rangle\subset\mathbb Q[x]$ and $I_{\mathbb{F}_p}:=\langle\overline f_1,\ldots,\overline f_m\rangle\subset\mathbb F_p[x]$.
Then
\[
 \HF_{\mathbb Q[x]/I}(d)
 \leq
 \HF_{\mathbb F_p[x]/I_{\mathbb{F}_p}}(d)
 \qquad(d\geq0).
\]
In particular, if $\HS_{\mathbb F_p[x]/I_{\mathbb{F}_p}}(z)=H_\tau^{\mathrm{Fr}}(z)$, then $\HS_{\mathbb Q[x]/I}(z)=H_\tau^{\mathrm{Fr}}(z)$.
\end{lemma}

\begin{proof}
For each $d\geq0$, form the integer matrix whose rows are the coefficient vectors, with respect to the degree-$d$ monomials, of all polynomials $u f_j$, where $d_j\leq d$ and $u$ ranges over the monomials of degree $d-d_j$.
Its rank over $\mathbb Q$ is at least the rank of its reduction modulo $p$, which gives the first assertion.

Suppose that $\HS_{\mathbb F_p[x]/I_{\mathbb{F}_p}}(z)=H_\tau^{\mathrm{Fr}}(z)$.
The first assertion shows that the coefficients of $\HS_{\mathbb Q[x]/I}(z)$ are no larger than those of $H_\tau^{\mathrm{Fr}}(z)$.
On the other hand, Fr\"oberg's inequality gives $\HS_{\mathbb Q[x]/I}(z)\succeq_{\mathrm{lex}}H_\tau^{\mathrm{Fr}}(z)$.
As in the proof of Proposition~\ref{prop:Frob-one-example}, these two inequalities imply $\HS_{\mathbb Q[x]/I}(z)=H_\tau^{\mathrm{Fr}}(z)$.
\end{proof}

Combining Lemma~\ref{lem:hilbert-series-lifting} with Theorem~\ref{thm:fixed-degreewise-specialization} gives the following criterion using an example over a prime field.

\begin{corollary}[Lifting from a prime field to characteristic zero]\label{cor:arithmetic-lifting}
Let $p$ be a prime, let $\tau=(n;d_1,\ldots,d_m)$, and let $\overline{\bm F}=(\overline f_1,\ldots,\overline f_m)$ with $\overline f_j\in\mathbb F_p[x]_{d_j}$.
Put $I_{\mathbb{F}_p}:=\langle\overline{\bm F}\rangle\subset\mathbb F_p[x]$.
Suppose that $\HS_{\mathbb F_p[x]/I_{\mathbb{F}_p}}(z)=H_\tau^{\mathrm{Fr}}(z)$ and that $\ini_\prec(I_{\mathbb{F}_p})$ is weakly reverse lexicographic.
For each $j$, choose any homogeneous polynomial $f_j\in\mathbb Z[x]_{d_j}$ whose reduction modulo $p$ is $\overline f_j$, and put $I:=\langle f_1,\ldots,f_m\rangle\subset\mathbb Q[x]$.
Then
\[
 \HS_{\mathbb Q[x]/I}(z)=H_\tau^{\mathrm{Fr}}(z).
\]
Moreover, $\GB_\prec(I)\subset\mathbb Z_{(p)}[x]$ and
\[
 \pi_p\bigl(\GB_\prec(I)\bigr)=\GB_\prec(I_{\mathbb{F}_p}),
 \qquad
 \LM_\prec\bigl(\GB_\prec(I)\bigr)
 =\LM_\prec\bigl(\GB_\prec(I_{\mathbb{F}_p})\bigr),
\]
where $\pi_p:\mathbb Z_{(p)}[x]\to\mathbb F_p[x]$ is coefficientwise reduction modulo $p$.
Consequently, Fr\"oberg's and Moreno-Soc\'ias' conjectures hold for the type $\tau$ over every field of characteristic zero.
\end{corollary}

\begin{proof}
By Lemma~\ref{lem:hilbert-series-lifting}, we have $\HS_{\mathbb Q[x]/I}(z)=H_\tau^{\mathrm{Fr}}(z)$.
Now apply Theorem~\ref{thm:fixed-degreewise-specialization} with $B=\mathbb Z_{(p)}$, $L=\mathbb Q$, $k=\mathbb F_p$, and $I_B=\langle f_1,\ldots,f_m\rangle\subset B[x]$.
The hypotheses hold in every degree.
Choosing $D\geq0$ at least as large as the degree of every element of the two reduced Gr\"obner bases gives the assertions about coefficients, specialization, and leading monomials.
In particular, $\ini_\prec(I)$ is weakly reverse lexicographic.
Both this property and the Hilbert series are preserved under extension of the coefficient field from $\mathbb Q$ to any field of characteristic zero.
The final assertion follows from Corollary~\ref{cor:Frob-MS}.
\end{proof}

\begin{remark}[Relation with modular Gr\"obner basis computation]
\label{rem:modular-computation}
Let $p$ be a prime, fix a homogeneous sequence with integer coefficients, and use the notation of Corollary~\ref{cor:arithmetic-lifting}.
Put $\mathcal G:=\GB_\prec(I)$ and $\mathcal H:=\GB_\prec(I_{\mathbb{F}_p})$.
Following Arnold~\cite[Definitions~5.1 and~5.2]{Arnold}, the prime $p$ is called \emph{Hilbert lucky} if the quotients $\mathbb Q[x]/I$ and $\mathbb F_p[x]/I_{\mathbb{F}_p}$ have the same Hilbert function, and \emph{lucky} if $\LM_\prec(\mathcal G)=\LM_\prec(\mathcal H)$.
The former condition does not imply the latter in general.

Suppose that $p$ is Hilbert lucky and that $\ini_\prec(I_{\mathbb{F}_p})$ is weakly reverse lexicographic.\footnote{Arnold assumes primitive integer generators.
Under the equality of the Hilbert functions, this causes no restriction here.
Replacing the nonzero generators by their primitive parts preserves $I$ and can only enlarge $I_{\mathbb{F}_p}$.
Lemma~\ref{lem:hilbert-series-lifting} then forces the resulting ideal over $\mathbb F_p$ to have the same Hilbert function as $I_{\mathbb{F}_p}$, so it equals $I_{\mathbb{F}_p}$.}
We claim that $p$ is lucky.
Suppose otherwise, and let $d$ be the smallest degree such that $\LM_\prec(\mathcal G_d)\ne\LM_\prec(\mathcal H_d)$.
The equality of the leading monomials in smaller degrees gives $B_d(I)=B_d(I_{\mathbb{F}_p})$.
Write this common set as $\{t_1,\ldots,t_w\}$, where $t_1\succ\cdots\succ t_w$.
By \eqref{eq:Bd-HF} and the equality of the Hilbert functions, we have $\#\mathcal G_d=\#\mathcal H_d=:v>0$.
Lemma~\ref{lem:wrl-initial-segment} gives $\LM_\prec(\mathcal H_d)=\{t_1,\ldots,t_v\}$.
Since $\LM_\prec(\mathcal G_d)\ne\{t_1,\ldots,t_v\}$ and both sets have $v$ elements, the smallest monomial $u$ in $\LM_\prec(\mathcal G_d)$ satisfies $u\prec t_v$.
However, $\prec$ is a graded order, and the leading monomials agree in all degrees smaller than $d$.
Thus \cite[Theorem~5.6]{Arnold}, applied with the leading monomials arranged in increasing order, implies that $t_v$ is not larger than $u$, a contradiction.
Hence $p$ is lucky.
It follows from \cite[Theorems~5.12 and~6.2]{Arnold} that $\mathcal G\subset\mathbb Z_{(p)}[x]$ and $\pi_p(\mathcal G)=\mathcal H$.

Thus the assertions about Hilbert series and reduced Gr\"obner bases in Corollary~\ref{cor:arithmetic-lifting} can also be deduced from Lemma~\ref{lem:hilbert-series-lifting} and Arnold's results.
In the proof of the corollary, we obtain the assertions about reduced Gr\"obner bases directly from Theorem~\ref{thm:fixed-degreewise-specialization}, which applies both to reduction modulo a prime and to specialization of coefficient parameters.
The hypotheses of the corollary can be verified entirely over $\mathbb F_p$, without first computing a Gr\"obner basis over $\mathbb Q$.
Recovering the rational coefficients still requires further computation.

Conversely, fix a type $\tau$ and suppose that both conjectures hold over $\mathbb Q$ for this type.
By Corollary~\ref{cor:Frob-MS}, there is a sequence over $\mathbb Q$ satisfying the simultaneous criterion.
After rescaling its nonzero entries, we may assume that they are primitive polynomials in $\mathbb Z[x]$.
By \cite[Theorem~5.13]{Arnold}, all but finitely many primes are lucky for this sequence.
For each such prime $p$, the initial ideals over $\mathbb Q$ and $\mathbb F_p$ have the same minimal monomial generators.
Hence the reduction of the sequence modulo $p$ also satisfies the simultaneous criterion.
After extension of scalars, Corollary~\ref{cor:Frob-MS} shows that both conjectures hold for $\tau$ over every infinite field of characteristic $p$.
\end{remark}

We illustrate Corollaries~\ref{cor:Frob-MS} and~\ref{cor:arithmetic-lifting} with the following computational results.
For $8\leq m\leq16$, consider the type $\tau_m=(8;\underbrace{2,\ldots,2}_{m})$.

\begin{proposition}[Computational verification]\label{prop:computed-types}
For every $m\in\{8,9,\ldots,16\}$, there exists a sequence $\bm F_m\in\mathbb F_{31}[x_1,\ldots,x_8]_2^m$ such that $\HS_{\mathbb F_{31}[x_1,\ldots,x_8]/I_{\bm F_m}}(z)=H_{\tau_m}^{\mathrm{Fr}}(z)$ and $\ini_\prec(I_{\bm F_m})$ is weakly reverse lexicographic.
Consequently, Fr\"oberg's and Moreno-Soc\'ias' conjectures hold over $\overline{\mathbb F}_{31}$ and over every field of characteristic zero for all nine types $\tau_m$.
\end{proposition}

\begin{proof}
For each $m$, the first assertion was checked in \textsc{Magma} \cite{BosmaCannonPlayoust} by direct computation of the Hilbert series and the reduced Gr\"obner basis.
Code for searching for and verifying such examples is given in Appendix~\ref{app:magma}.
The properties are unchanged after extending scalars from $\mathbb F_{31}$ to $\overline{\mathbb F}_{31}$.
The assertion over $\overline{\mathbb F}_{31}$ follows from Corollary~\ref{cor:Frob-MS}.
The characteristic-zero assertion follows from Corollary~\ref{cor:arithmetic-lifting}.
\end{proof}

\subsection{Reduction to complete intersections}\label{sec:ci}

Fix $m\geq1$ and positive integers $d_1,\ldots,d_m$.
Pardue's argument~\cite[Section~4, p.~583]{Pardue} shows that if Moreno-Soc\'ias' conjecture holds for the type $(m;d_1,\ldots,d_m)$, then it holds for the type $(n;d_1,\ldots,d_m)$ for every $n\geq1$.
Combining this reduction with Theorem~\ref{thm:generic-HS-criterion} gives the following complete-intersection criterion.

\begin{corollary}[Reduction to one complete intersection]\label{cor:CI-construction}
Fix $m\geq1$ and positive integers $d_1,\ldots,d_m$, and put $S:=K[x_1,\ldots,x_m]$.
Suppose that there exists a homogeneous regular sequence
$\bm F=(f_1,\ldots,f_m)$ in $S$, with $\deg f_i=d_i$ for each $i$, such that $\ini_\prec(I_{\bm F})$ is weakly reverse lexicographic.
Then Moreno-Soc\'ias' conjecture holds for the type $(n;d_1,\ldots,d_m)$ for every $n\geq1$.
In particular, if such a regular sequence exists for every $m\geq1$ and every choice of positive degrees $d_1,\ldots,d_m$, then Moreno-Soc\'ias' conjecture holds for every type.
\end{corollary}

\begin{proof}
Put $\tau=(m;d_1,\ldots,d_m)$.
Both $\bm F$ and the universal generic sequence $\bm G_\tau$ are regular sequences of degrees $d_1,\ldots,d_m$.
Hence we have $\HS_{S/I_{\bm F}}(z)=\frac{\prod_{i=1}^m(1-z^{d_i})}{(1-z)^m} =H_\tau^{\mathrm{gen}}(z)$; see, for example,
\cite[Chapter~17]{Eisenbud}.
Since $\ini_\prec(I_{\bm F})$ is weakly reverse lexicographic by hypothesis, $\bm F$ satisfies condition~\textup{(iii)} of
Theorem~\ref{thm:generic-HS-criterion}.
Thus Moreno-Soc\'ias' conjecture holds for $\tau$.
The assertion for every $n\geq1$ now follows from Pardue's reduction above.
\end{proof}

For example, the case $m=8$ of Proposition~\ref{prop:computed-types} gives a complete intersection of $8$ quadrics in $8$ variables.
By Corollary~\ref{cor:CI-construction}, Moreno-Soc\'ias' conjecture therefore holds for $8$ quadrics in any number of variables over every field of characteristic zero and every infinite field of characteristic $31$.

\begin{problem}[Uniform complete-intersection construction]\label{prob:CI}
Given $n$ and $d_1,\ldots,d_n$, construct explicitly one regular sequence of these degrees whose initial ideal with respect to a graded reverse lexicographic order is weakly reverse lexicographic.
\end{problem}

The existence of a weakly reverse lexicographic monomial ideal with the Hilbert function of such a complete intersection follows from \cite[Corollary~6]{Pardue}; for an explicit construction when all $d_i\geq 2$, see \cite[Theorem~4.1]{BertoneCioffi}.
Problem~\ref{prob:CI} asks for a regular sequence of the prescribed degrees whose initial ideal is this monomial ideal.

Suppose that $K$ has characteristic zero.
One way to seek examples for Problem~\ref{prob:CI} is to consider sequences $(L_1^{d_1},\ldots,L_n^{d_n})$, where $L_1,\ldots,L_n$ are linearly independent linear forms in $R$.
Every such sequence is regular, since it is obtained from
$(x_1^{d_1},\ldots,x_n^{d_n})$ by an invertible linear change of coordinates.
Thus, writing $I_0=\langle x_1^{d_1},\ldots,x_n^{d_n}\rangle$, we ask whether $\ini_\prec(g(I_0))$ is weakly reverse lexicographic for some $g\in\mathrm{GL}_n(K)$.
Cho and Park~\cite{ChoPark} proved that $\gin_\prec(I_0)$ is weakly reverse lexicographic (almost reverse lexicographic in their terminology) if the degrees are indexed so that
\begin{equation}\label{eq:ChoPark-condition}
 d_i>\sum_{j=1}^{i-1}d_j-i+1
 \qquad(4\leq i\leq n).
\end{equation}
Since $\ini_\prec(g(I_0))=\gin_\prec(I_0)$ for all $g$ in a nonempty Zariski-open subset of $\mathrm{GL}_n(K)$, any such $g$ gives a regular sequence with the required initial ideal under \eqref{eq:ChoPark-condition}.

Cho and Park also deduced Moreno-Soc\'ias' conjecture for generic complete intersections of degrees satisfying \eqref{eq:ChoPark-condition}.
Capaverde and Gao~\cite{CapaverdeGao} proved the conjecture for generic complete intersections under the bound $d_i\geq\sum_{j=1}^{i-1}d_j-i-1$ for $2\leq i\leq n$, subject to the ordering assumptions in their theorem.
Their result concerns generic sequences in the full coefficient space $R_{d_1}\times\cdots\times R_{d_n}$; it does not assert that the same initial ideal can be obtained from powers of linearly independent linear forms.

In the case of $n=4$, results of Harima, Sakaki, and Wachi~\cite{HarimaSakakiWachi} illustrate both the scope and the limitation of using such powers.
They proved that $\gin_\prec(I_0)$ is weakly reverse lexicographic if at least one of the four degrees $d_1,d_2,d_3,d_4$ is $2$.
This includes the degree sequence $(2,2,2,2)$, which does not satisfy \eqref{eq:ChoPark-condition}.
On the other hand, their results imply that if $(d_1,d_2,d_3,d_4)=(3,3,3,3)$, then $\ini_\prec(g(I_0))$ is not weakly reverse lexicographic for any $g\in\mathrm{GL}_4(K)$.
Thus a construction that solves Problem~\ref{prob:CI} for all degree sequences cannot be restricted to powers of linearly independent linear forms.

\section{Sequences with a prescribed projective zero}\label{sec:dimone}

We now consider homogeneous polynomial sequences with a prescribed common projective zero.
We first apply the specialization theorem of Section~\ref{sec:specialization} to the corresponding coefficient space without any restriction on the number of generators.
We then obtain a Krull-dimension-one analogue of the one-example criterion under the assumption $m\geq n-1$.

\subsection{The coefficient space $Z_{\bm p}$}

We return to the setting of Section~\ref{subsec:coeff_sp}.
Thus $K$ is an infinite field of arbitrary characteristic, and $\tau=(n;d_1,\ldots,d_m)$ is fixed.
Fix a point $\bm p=(p_1:\cdots:p_n)\in\mathbb P^{n-1}(K)$, where $\mathbb{P}^{n-1}(K)$ denotes the projective $(n-1)$-space over $K$.
We choose a nonzero vector $(p_1,\ldots,p_n)\in K^n$ representing $\bm p$.
For $d\geq1$, put $R_d^{\bm p}:=\{f\in R_d:f(\bm p)=0\}$, and define the linear subspace
\[
 Z_{\tau,\bm p}:=
 R_{d_1}^{\bm p}\times\cdots\times R_{d_m}^{\bm p}\subset V_\tau.
\]
Since $\tau$ is fixed, we write $Z_{\bm p}$ for $Z_{\tau,\bm p}$.
Thus $Z_{\bm p}$ consists of the sequences $\bm F=(f_1,\ldots,f_m)\in V_{\tau}$ satisfying $f_j(\bm p)=0$ for all $j$.
Recall that $A_\tau=K[\mathcal C_\tau]$ is the coefficient ring of the universal generic sequence on $V_\tau$.
For $1\leq j\leq m$, put $\ell_{j,\bm p}:= \sum_{t\in\mathbb T_{d_j}}t(\bm p)c_{j,t} \in A_\tau$, and define $A_\tau^{\bm p}:= A_\tau/\langle \ell_{1,\bm p},\ldots,\ell_{m,\bm p}\rangle$.
Choose and fix an index $\nu$ with $p_\nu\neq0$.
For each $j$, put $t_j:=x_\nu^{d_j}$.
Since the coefficient of $c_{j,t_j}$ in $\ell_{j,\bm p}$ is $p_\nu^{d_j}\neq0$, we can solve $\ell_{j,\bm p}=0$ for $c_{j,t_j}$, obtaining
\[
A_\tau^{\bm p} \cong K[c_{j,t}:1\leq j\leq m,\ t\in\mathbb T_{d_j}\smallsetminus\{x_\nu^{d_j}\}].
\]
In particular, $Z_{\bm p}$ is an affine space of codimension $m$ in $V_\tau$.

\begin{lemma}\label{lem:Zp-coordinate-ring}
With notation as above, let $I_K(Z_{\bm p})$ denote the ideal of all polynomials in $A_\tau$ that vanish at every point of $Z_{\bm p}$. Then $I_K(Z_{\bm p})  =\langle\ell_{1,\bm p},\ldots,\ell_{m,\bm p}\rangle$.
In particular, $A_\tau^{\bm p}$ is the coordinate ring of $Z_{\bm p}$.
\end{lemma}

\begin{proof}
Put $J:=\langle\ell_{1,\bm p},\ldots,\ell_{m,\bm p}\rangle$.
For every $\bm F=(f_1,\ldots,f_m)\in Z_{\bm p}$, we have $\pi_{\bm F}(\ell_{j,\bm p})=f_j(\bm p)=0$ for all $j$.
Thus $J\subseteq I_K(Z_{\bm p})$.
Conversely, let $h\in I_K(Z_{\bm p})$.
Under the above identification of $A_\tau/J$ with the polynomial ring in the remaining $N-m$ coefficient variables, let $\overline h$ be the polynomial corresponding to $h+J$.
For each $\bm u\in K^{N-m}$, assigning the remaining coefficients the values $\bm u$ and solving for the eliminated coefficients gives a point $\bm F_{\bm u}\in Z_{\bm p}$.
By construction, $\overline h(\bm u)=\pi_{\bm F_{\bm u}}(h)=0$.
Since $K$ is infinite, $\overline h$ is the zero polynomial.
Hence $h+J=0$ in $A_\tau/J$, so $h\in J$.
\end{proof}

For each $\bm F\in Z_{\bm p}$, the coefficient specialization
$\pi_{\bm F}:A_\tau\to K$ therefore factors through $A_\tau^{\bm p}$;
we denote the induced map $A_\tau^{\bm p}\to K$ by the same symbol.
Writing $\overline c_{j,t}$ for the image of $c_{j,t}$ in $A_\tau^{\bm p}$, we define
\[
g_j^{\bm p}:= \sum_{t\in\mathbb T_{d_j}}\overline c_{j,t}t\in A_{\tau}^{\bm{p}}[x_1,\ldots,x_n],
\qquad \bm G_\tau^{\bm p}:= (g_1^{\bm p},\ldots,g_m^{\bm p}).
\]
Then $g_j^{\bm p}(\bm p)=0$ for every $1\leq j \leq m$.
The sequence $\bm G_\tau^{\bm p}$ over $A_\tau^{\bm p}$ is the universal family of sequences on $Z_{\bm p}$ in the following sense:
For every commutative $K$-algebra $C$ and every sequence $\bm H=(h_1,\ldots,h_m)\in\prod_{j=1}^m C[x]_{d_j}$ satisfying $h_j(\bm p)=0$ for all $j$, there is a unique $K$-algebra homomorphism $\pi_{\bm H}:A_\tau^{\bm p}\to C$ such that $\pi_{\bm H}(g_j^{\bm p})=h_j$ for all $j$, where $\pi_{\bm H}$ is applied coefficientwise.
We use the same notation for the sequence over $L_\tau^{\bm p}:=\Frac(A_\tau^{\bm p})$ and call it the {\it universal generic sequence on $Z_{\bm p}$}.
We put
\[
 I_{\tau,\bm p}^{\mathrm{gen}} := \langle\bm G_\tau^{\bm p}\rangle \subset L_\tau^{\bm p}[x].
\]
Since $\bm p$ is a common projective zero of $\bm G_\tau^{\bm p}$, the quotient $L_\tau^{\bm p}[x]/I_{\tau,\bm p}^{\mathrm{gen}}$ has positive Krull dimension.
The same holds for $R/I_{\bm F}$ for every $\bm F\in Z_{\bm p}$.
The {\it generic Hilbert series on $Z_{\bm p}$} is defined by
\[
 H_{\tau,\bm p}^{\mathrm{gen}}(z)
 :=
 \HS_{L_\tau^{\bm p}[x]/I_{\tau,\bm p}^{\mathrm{gen}}}(z).
\]
For $\bm F\in Z_{\bm p}$, let $\mathcal H_{\tau,\bm p}^{\mathrm{gen}}(\bm{F})$ denote the property $\HS_{R/I_{\bm F}}(z)=H_{\tau,\bm p}^{\mathrm{gen}}(z)$.
We also call $\HF_{L_\tau^{\bm p}[x]/I_{\tau,\bm p}^{\mathrm{gen}}}$ the {\it generic Hilbert function on $Z_{\bm p}$}.
For every $\bm F\in Z_{\bm p}$ and every $d\geq0$, we have
\begin{equation}\label{eq:Zp-HF-bounds}
 1\leq \HF_{L_\tau^{\bm p}[x]/I_{\tau,\bm p}^{\mathrm{gen}}}(d)\leq \HF_{R/I_{\bm F}}(d).
\end{equation}
The lower bound follows from $\dim L_\tau^{\bm p}[x]/I_{\tau,\bm p}^{\mathrm{gen}} \geq 1$, whereas the upper bound follows from \cite[Appendix~B, Lemma~3]{Diem} applied to the specialization $\pi_{\bm F}:A_\tau^{\bm p}\to K$.

The following theorem is the analogue of Theorem~\ref{thm:generic-HS-criterion} for $Z_{\bm p}$.
It does not require an explicit formula for $H_{\tau,\bm p}^{\mathrm{gen}}(z)$.

\begin{theorem}[One-example criterion for a weakly reverse lexicographic initial ideal on $Z_{\bm p}$]\label{thm:Zp-generic-HS}
For the fixed type $\tau$ and point $\bm p$, the following are equivalent:
\begin{enumerate}[label=\textup{(\roman*)}]
 \item $\ini_\prec(I_{\tau,\bm p}^{\mathrm{gen}})$ is weakly reverse lexicographic;
 \item $\mathcal H_{\tau,\bm p}^{\mathrm{gen}}(\bm{F})\wedge\mathcal W({\bm F})$ holds for every $\bm{F}$ in a nonempty Zariski-open subset of $Z_{\bm p}$;
 \item there exists one $\bm F\in Z_{\bm p}$ satisfying $\mathcal H_{\tau,\bm p}^{\mathrm{gen}}(\bm{F})\wedge\mathcal W({\bm F})$.
\end{enumerate}
\end{theorem}

\begin{proof}
Since $A_\tau^{\bm p}$ is a polynomial ring over $K$, the proof of Proposition~\ref{prop:generic-stability} applies with $A_\tau$ and $\bm G_\tau$ replaced by $A_\tau^{\bm p}$ and $\bm G_\tau^{\bm p}$, respectively.
Hence there is a nonempty relatively Zariski-open subset of $Z_{\bm p}$ on which the specialized initial ideals have the same monomials as $\ini_\prec(I_{\tau,\bm p}^{\mathrm{gen}})$.
Thus (i) implies (ii), and (ii) implies (iii) immediately.

Assume (iii), and choose $\bm F\in Z_{\bm p}$ satisfying $\mathcal H_{\tau,\bm p}^{\mathrm{gen}}(\bm F)\wedge\mathcal W(\bm F)$.
Since $\pi_{\bm F}:A_\tau^{\bm p}\to K$ restricts to the identity on $K$, it is surjective.
Thus $\mathfrak n:=\ker(\pi_{\bm F})$ is a maximal ideal of $A_\tau^{\bm p}$ with $A_\tau^{\bm p}/\mathfrak n\cong K$.
Put $B:=(A_\tau^{\bm p})_{\mathfrak n}$, and let $I_B:=\langle g_1^{\bm p},\ldots,g_m^{\bm p}\rangle\subset B[x]$.
Then $B$ is a Noetherian local domain with fraction field $L_\tau^{\bm p}$ and residue field $K$.
The ideals $I_L$ and $I_k$ of Theorem~\ref{thm:fixed-degreewise-specialization} are $\langle\bm G_\tau^{\bm p}\rangle\subset L_\tau^{\bm p}[x]$ and $I_{\bm F}\subset R$, respectively.
The condition $\mathcal H_{\tau,\bm p}^{\mathrm{gen}}(\bm{F})$ gives the required equality of Hilbert functions in every degree, while $\mathcal W({\bm F})$ gives the weakly reverse lexicographic hypothesis.
Choose $D\geq0$ at least as large as the degree of every element of the two reduced Gr\"obner bases.
Theorem~\ref{thm:fixed-degreewise-specialization} then shows that the two reduced Gr\"obner bases have the same leading monomials.
Hence the initial ideal of the generic sequence is weakly reverse lexicographic, proving (i).
\end{proof}

\subsection{An explicit Hilbert-series criterion}

The following lemma explains the condition $m\geq n-1$
in the Krull-dimension-one setting.

\begin{lemma}\label{lem:Zp-dimension-one}
There exists a sequence $\bm F\in Z_{\bm p}$ with
$\dim R/I_{\bm F}=1$ if and only if $m\geq n-1$.
Moreover, if $m\geq n-1$, then
$\dim\bigl(L_\tau^{\bm p}[x]/I_{\tau,\bm p}^{\mathrm{gen}}\bigr)=1$.
\end{lemma}

\begin{proof}
If $\bm F\in Z_{\bm p}$ satisfies $\dim R/I_{\bm F}=1$, then Krull's height theorem gives $1= \dim R/I_{\bm F}\geq n-m$, so $m\geq n-1$.
Conversely, assume that $m\geq n-1$.
Choose a $K$-linear basis $\ell_1,\ldots,\ell_{n-1}$ of $R_1^{\bm p}$ and a linear form $\ell_n\in R_1$ with $\ell_n(\bm p)=1$.
Then $\ell_1,\ldots,\ell_n$ are algebraically independent over $K$, and $R=K[\ell_1,\ldots,\ell_n]$.
Thus, for $P:=\langle\ell_1,\ldots,\ell_{n-1}\rangle$, we have $R/P\cong K[\ell_n]$, so that $P$ is a prime homogeneous ideal with $\dim R/P=1$.
It therefore suffices to construct $\bm F\in Z_{\bm p}$ such that $\sqrt{I_{\bm F}}=P$.
Set $f_j:=\ell_j^{d_j}$ for $1\leq j\leq n-1$, and choose any $f_j\in P\cap R_{d_j}$ for $n\leq j\leq m$.
The resulting sequence $\bm F=(f_1,\ldots,f_{m})$ belongs to $Z_{\bm p}$.
Then $I_{\bm F}\subseteq P$, and so $\sqrt{I_{\bm F}}\subseteq P$ as $P$ is prime.
On the other hand, since $\ell_j^{d_j}\in I_{\bm F}$ for $1\leq j\leq n-1$, we have $\sqrt{I_{\bm F}}\supseteq P$.
Thus $\sqrt{I_{\bm F}}=P$.
Since passing to the radical does not change the Krull dimension of the quotient, we have $\dim R/I_{\bm F}=\dim R/P=1$.

Since $\dim R/I_{\bm F}=1$, the Hilbert function of $R/I_{\bm F}$ is eventually equal to a positive constant.
By \eqref{eq:Zp-HF-bounds}, the Hilbert polynomial of $L_\tau^{\bm p}[x]/I_{\tau,\bm p}^{\mathrm{gen}}$ is therefore bounded above and below by positive constants for all sufficiently large integers.
It is thus a positive constant, so the generic quotient has Krull dimension one.
\end{proof}

We henceforth assume that $m\geq n-1$, and next introduce a predicted Hilbert series for sequences with a prescribed projective zero and establish a one-example criterion for its validity.
For two power series with real coefficients, $\max\{A(z),B(z)\}$ denotes coefficientwise maximum.
Put
\begin{equation}\label{eq:one-dimensional-series}
H_\tau^{\mathrm{Fr},(1)}(z):= \max\left\{ H_\tau^{\mathrm{Fr}}(z), \frac{1}{1-z} \right\}.
\end{equation}
For $m=n-1$, the coefficients of $H_\tau^{\mathrm{Fr},(1)}(z)$ are eventually $\prod_{j=1}^{n-1}d_j$, whereas for $m\geq n$ they are eventually $1$.
Thus any homogeneous quotient with this Hilbert series has Krull dimension one.

\begin{proposition}[One-example criterion for the predicted Hilbert series on $Z_{\bm p}$]
\label{prop:H1-one-example}
For a fixed type $\tau=(n;d_1,\ldots,d_m)$ with $m\geq n-1$, the following are equivalent:
\begin{enumerate}[label=\textup{(\roman*)}]
 \item $H_{\tau,\bm p}^{\mathrm{gen}}(z) =H_\tau^{\mathrm{Fr},(1)}(z)$;
 \item the equality $\HS_{R/I_{\bm F}}(z)=H_\tau^{\mathrm{Fr},(1)}(z)$ holds for every $\bm{F}$ in a nonempty Zariski-open subset of $Z_{\bm p}$;
 \item there exists one $\bm F\in Z_{\bm p}$ such that $\HS_{R/I_{\bm F}}(z)=H_\tau^{\mathrm{Fr},(1)}(z)$.
\end{enumerate}
\end{proposition}

\begin{proof}
The implication (i) $\Rightarrow$ (ii) follows from the generic stability on $Z_{\bm p}$ established in the proof of Theorem~\ref{thm:Zp-generic-HS}, while (ii) $\Rightarrow$ (iii) is immediate.

Assume (iii), and choose such a sequence $\bm F$.
By \eqref{eq:Zp-HF-bounds}, we have $1\leq \HF_{L_\tau^{\bm p}[x]/I_{\tau,\bm p}^{\mathrm{gen}}}(d)\leq\HF_{R/I_{\bm F}}(d)$ for every $d\geq0$.
Fr\"oberg's inequality \cite{Froberg} gives $H_{\tau,\bm p}^{\mathrm{gen}}(z)\succeq_{\mathrm{lex}}H_\tau^{\mathrm{Fr}}(z)$.
Suppose that $\HS_{L_\tau^{\bm p}[x]/ \langle\bm G_\tau^{\bm p}\rangle}(z)$ differs from $\HS_{R/I_{\bm F}}(z)=H_\tau^{\mathrm{Fr},(1)}(z)$, and let $d$ be the first degree at which their coefficients differ.
Since $\HF_{L_\tau^{\bm p}[x]/\langle\bm G_\tau^{\bm p}\rangle}(d)\leq\HF_{R/I_{\bm F}}(d)$, the inequality must be strict in degree $d$.
As $\HF_{L_\tau^{\bm p}[x]/ \langle\bm G_\tau^{\bm p}\rangle}(d) \geq 1$, the coefficient of $z^d$ in $H_\tau^{\mathrm{Fr},(1)}(z)$ must be greater than $1$.
Thus, the definition of Fr\"oberg's truncated series~\eqref{eq:Frob-series} implies that $H_\tau^{\mathrm{Fr},(1)}(z)$ and $H_\tau^{\mathrm{Fr}}(z)$ agree through degree $d$.
Hence the generic Hilbert series $\HS_{L_\tau^{\bm p}[x]/ \langle\bm G_\tau^{\bm p}\rangle}(z)$ agrees with $H_\tau^{\mathrm{Fr}}(z)$ in every degree below $d$ but has a smaller coefficient in degree $d$.
This contradicts $H_{\tau,\bm p}^{\mathrm{gen}}(z)\succeq_{\mathrm{lex}}H_\tau^{\mathrm{Fr}}(z)$.
Hence $H_{\tau,\bm p}^{\mathrm{gen}}(z)=H_\tau^{\mathrm{Fr},(1)}(z)$, proving (i).
\end{proof}

\subsection{A one-example criterion in Krull dimension one}

For $\bm F\in Z_{\bm p}$, let $\mathcal H_{\tau}^{(1)}(\bm F)$ denote the property
\begin{equation*}
 \HS_{R/I_{\bm F}}(z)=H_\tau^{\mathrm{Fr},(1)}(z).
\end{equation*}
Combining Proposition~\ref{prop:H1-one-example} with Theorem~\ref{thm:Zp-generic-HS} gives the following criterion.

\begin{theorem}[Simultaneous one-example criterion on $Z_{\bm p}$ in Krull dimension one]\label{thm:dimone-criterion}
For every fixed type $\tau=(n;d_1,\ldots,d_m)$ with $m\geq n-1$, the following are equivalent:
\begin{enumerate}[label=\textup{(\roman*)}]
 \item $H_{\tau,\bm p}^{\mathrm{gen}}(z) =H_\tau^{\mathrm{Fr},(1)}(z)$ and $\ini_\prec(I_{\tau,\bm p}^{\mathrm{gen}})$ is weakly reverse lexicographic;
 \item $\mathcal H_{\tau}^{(1)}(\bm F)\wedge\mathcal W(\bm F)$ holds for every $\bm F$ in a nonempty Zariski-open subset of $Z_{\bm p}$;
 \item there exists one $\bm F\in Z_{\bm p}$ satisfying
 $\mathcal H_{\tau}^{(1)}(\bm F)\wedge\mathcal W(\bm F)$.
\end{enumerate}
\end{theorem}

\begin{proof}
The implication (i) $\Rightarrow$ (ii) follows from the generic stability on $Z_{\bm p}$ established in the proof of Theorem~\ref{thm:Zp-generic-HS}, while (ii) $\Rightarrow$ (iii) is immediate.

Assume (iii), and choose such a sequence $\bm F$.
Proposition~\ref{prop:H1-one-example} gives $H_{\tau,\bm p}^{\mathrm{gen}}(z) =H_\tau^{\mathrm{Fr},(1)}(z)$.
Thus $\bm F$ satisfies $\mathcal H_{\tau,\bm p}^{\mathrm{gen}}(\bm F)$.
Theorem~\ref{thm:Zp-generic-HS} then shows that $\ini_\prec(I_{\tau,\bm p}^{\mathrm{gen}})$ is weakly reverse lexicographic, proving (i).
\end{proof}

The preceding criteria hold for every $\bm p=(p_1:\cdots:p_n)\in\mathbb P^{n-1}(K)$, without any restriction on its coordinates.
The criteria do not require a change of coordinates taking $\bm p$ to $\bm o=(0:\cdots:0:1)$, nor do they require $K$ to be algebraically closed.

For the conjecture below, however, the position of $\bm p$ relative to the fixed graded reverse lexicographic order matters.
Let $\bm o=(0:\cdots:0:1)$.
If $p_n\neq0$, we may normalize $p_n=1$.
The substitution $x_i\mapsto x_i-p_ix_n$ for $1\leq i<n$, with $x_n$ fixed, induces an isomorphism $Z_{\bm o}\to Z_{\bm p}$.
It preserves leading monomials, since each monomial is mapped to itself plus monomials with larger $x_n$-exponents, which are smaller in the graded reverse lexicographic order.
Being a graded automorphism, it also preserves Hilbert series.
Consequently, whenever $p_n\neq0$, the equivalent conditions of Theorem~\ref{thm:dimone-criterion} hold for $\bm p$ if and only if they hold for $\bm o$.

\begin{conjecture}[Krull-dimension-one variant]\label{conj:dimone}
For every type $\tau=(n;d_1,\ldots,d_m)$ with $m\geq n-1$, the equivalent conditions of Theorem~\ref{thm:dimone-criterion} hold for $\bm p=\bm o$.
Equivalently, they hold for every $\bm p\in\mathbb P^{n-1}(K)$ with $p_n\neq0$.
\end{conjecture}

\begin{remark}
The restriction $p_n\neq0$ in Conjecture~\ref{conj:dimone} cannot be omitted in general.
For example, take $\tau=(2;1)$ and $\bm p=(1:0)$.
Every linear form vanishing at $\bm p$ is a scalar multiple of $x_2$, so $I_{\tau,\bm p}^{\mathrm{gen}}=\langle x_2\rangle$.
This ideal is not weakly reverse lexicographic, since $x_1\succ x_2$ but $x_1\notin\langle x_2\rangle$.
For $\bm o=(0:1)$, by contrast, the corresponding ideal is $\langle x_1\rangle$, which is weakly reverse lexicographic.
In both cases the quotient has the predicted Hilbert series $1/(1-z)$.
\end{remark}

\section*{Acknowledgements}
The authors thank the editor of the previous submission for suggesting that the original manuscript be divided into two focused papers.
This work was supported by JSPS Grant-in-Aid for Young Scientists 23K12949, JSPS Grant-in-Aid for Scientific Research (C) 21K03377, (B) 25K00140, and 26K06753, and JST CREST Grant Number JPMJCR2113.

\section*{Data availability}
The \textsc{Magma} code used to search for the examples in Proposition~\ref{prop:computed-types} is included in Appendix~\ref{app:magma}.
The coefficient matrices, the corresponding polynomials in expanded form, the output of the original search, and a script for verifying all nine examples are available at \url{https://github.com/Momonari-Kudo/moreno-socias-magma-examples}.

\appendix

\section{Specialization via marked families}\label{app:marked}

The purpose of this appendix is to give an alternative proof of the specialization results from a module-theoretic viewpoint.
Unless otherwise specified, $\prec$ denotes the graded reverse lexicographic order with $x_1\succ\cdots\succ x_n$, as in Section~\ref{sec:prelim}.
Our ordering of the variables is the reverse of that used in~\cite{LR16}.
The proof is based on the following two observations.

\begin{enumerate}[label=\textup{(\roman*)}]
\item Equality of the Hilbert functions of the generic and special fibers is equivalent to membership in the marked family associated with the initial ideal of the special fiber.
This equivalence requires no stability assumption on the initial ideal.
The nontrivial implication is proved in Lemma~\ref{lem:graded-local-free}, and the equivalence is summarized in Theorem~\ref{thm:marked-specialization}.

\item If the initial ideal of the special fiber is weakly reverse lexicographic, the corresponding unique marked basis is a reduced Gr\"obner basis (Lemma~\ref{lem:wrl-marked-gb}).
\end{enumerate}
By the base-change property of marked bases (Proposition~\ref{prop:marked-base-change}), the reduced Gr\"obner bases of the generic and special fibers are then obtained from the same basis over the local coefficient ring, with identical leading monomials (Theorem~\ref{thm:marked-specialization}).

We first recall the necessary results on marked bases and base change, following Lella and Roggero~\cite{LR16}.
We then prove the specialization theorem and derive its finite-degree and prescribed-zero variants.

\subsection{Marked bases and marked families}

Marked bases generalize reduced Gr\"obner bases by prescribing head monomials without requiring them to be leading monomials with respect to a monomial order.
Their defining property is that the monomials outside a fixed monomial ideal form a free basis of the quotient over the coefficient ring.
We recall the definitions and results needed below, following Lella and Roggero~\cite{LR16}.

Let $C$ be a commutative Noetherian ring, and let $C[x]=C[x_1,\ldots,x_n]$.
For a monomial ideal $J$, denote by $\mathcal B_J$ its minimal monomial generating set and put
\[
\mathcal N(J):=\{x^\beta\mid x^\beta\notin J\},
\qquad
\mathcal N(J)_d:=\mathcal N(J)\cap\mathbb T_d,
\qquad
\langle\mathcal N(J)\rangle_C
:=
\bigoplus_{x^\beta\in\mathcal N(J)}Cx^\beta.
\]
We regard $J$ as a monomial ideal determined by its generating monomials, independently of the coefficient ring, and use the same notation for its extension to $C[x]$.

Throughout this appendix, by a \emph{reduced Gr\"obner basis} over a coefficient ring, we mean a monic Gr\"obner basis whose leading monomials form the minimal monomial generating set of the initial ideal and whose nonleading monomials lie outside that ideal.

A \emph{marked polynomial} is a polynomial with a specified monomial in its support, called its \emph{head term} and denoted by $\operatorname{Ht}(f)$.
It is called \emph{monic} if the coefficient of its head term is $1$.
There need not exist a monomial order for which the head term is the leading monomial.

\begin{definition}[$J$-marked set and $J$-marked basis]
A \emph{$J$-marked set} is a collection $\mathcal F_J=\{f_\alpha\mid x^\alpha\in\mathcal B_J\}$ of monic marked polynomials of the form
\begin{equation}\label{eq:marked-poly}
f_\alpha
=
x^\alpha-
\sum_{x^\beta\in\mathcal N(J)_{|\alpha|}}
c_{\alpha\beta}x^\beta
\qquad(c_{\alpha\beta}\in C),
\end{equation}
where $\operatorname{Ht}(f_\alpha)=x^\alpha\in \mathcal{B}_J$.
The set $\mathcal{F}_J$ is called a \emph{$J$-marked basis} if, for $I:=\langle\mathcal F_J\rangle\subset C[x]$, one has
\[
C[x]=I\oplus_C\langle\mathcal N(J)\rangle_C.
\]
Equivalently, the residue classes of $\mathcal N(J)$ form a free $C$-basis of $C[x]/I$.
\end{definition}

\begin{definition}[$J$-marked family]
For a monomial ideal $J$, put
\[
\mathbf{Mf}_J(C)
:=
\left\{
I\subset C[x]\text{ homogeneous ideal}
\;\middle|\;
C[x]=I\oplus_C\langle\mathcal N(J)\rangle_C
\right\}.
\]
We call $\mathbf{Mf}_J(C)$ the \emph{$J$-marked family over $C$}.
\end{definition}

Recall that a monomial ideal $J$ is \emph{strongly stable} if, for every $x^\alpha\in J$ and every $i<j$ with $x_j\mid x^\alpha$, one has $x_i(x^\alpha/x_j)\in J$.
Every weakly reverse lexicographic monomial ideal is strongly stable, but the converse need not hold in general:
The ideal $J=\langle x_1^2,x_1x_2,x_1x_3\rangle\subset K[x_1,x_2,x_3]$ is strongly stable, but $x_1x_3$ is a minimal generator of $J$, whereas the larger monomial $x_2^2$ does not belong to $J$.

The following characterization is known.

\begin{proposition}[{\cite[Proposition~2.1]{LR16}}]\label{prop:marked-basis-I}
Let $J$ be a strongly stable monomial ideal and let $I\subset C[x]$ be a homogeneous ideal.
Then $I\in\mathbf{Mf}_J(C)$ if and only if $I$ is generated by a $J$-marked basis.
Moreover, such a $J$-marked basis is uniquely determined by $I$.
\end{proposition}

The ``if'' direction follows from the definition of a $J$-marked basis.
For the converse, the point is to show that the marked polynomials determined by the direct sum decomposition generate $I$.
The construction is compatible with changes of coefficient rings.

\begin{proposition}[Base change of marked families~{\cite[Proposition~2.3]{LR16}}]\label{prop:marked-base-change}
Let $J$ be a strongly stable monomial ideal, let $C,C'$ be commutative Noetherian rings, and let $\varphi:C\to C'$ be a ring homomorphism.
We also denote by $\varphi:C[x]\to C'[x]$ its coefficientwise extension, which fixes the variables.
Suppose that $I\in\mathbf{Mf}_J(C)$, and let $I_{C'}:=IC'[x]$ denote the ideal of $C'[x]$ generated by $\varphi(I)$.

Then $I_{C'}\in\mathbf{Mf}_J(C')$, and its unique $J$-marked basis is obtained by applying $\varphi$ coefficientwise to the unique $J$-marked basis of $I$.
\end{proposition}

In particular, no flatness assumption on $C'$ over $C$ is required.
Indeed, the defining decomposition $C[x]=I\oplus_C\langle\mathcal N(J)\rangle_C$ is preserved by tensoring with $C'$.

We will use the following consequence for generic and special fibers.

\begin{corollary}[Generic and special fibers]\label{cor:marked-BLk}
Let $(B,\mathfrak m)$ be a Noetherian local domain, and put $L:=\operatorname{Frac}(B)$ and $k:=B/\mathfrak m$.
Let $J$ be a strongly stable monomial ideal, and suppose that $I_B\in\mathbf{Mf}_J(B)$.
Set $I_L:=I_BL[x]$ and $I_k:=I_Bk[x]$.

Then $I_L\in\mathbf{Mf}_J(L)$ and $I_k\in\mathbf{Mf}_J(k)$.
Moreover, their unique $J$-marked bases are obtained from the unique $J$-marked basis of $I_B$ by applying the embedding $B\to L$ and 
the projection $B\to k$ coefficientwise, respectively.
In particular, the coefficients $c_{\alpha\beta}\in B$ in \eqref{eq:marked-poly} become $c_{\alpha\beta}\in L$ and $\overline{c_{\alpha\beta}}\in k$, respectively.
\end{corollary}

\begin{proof}
Apply Proposition~\ref{prop:marked-base-change} to the natural homomorphisms $B\to L$ and $B\to k$.
\end{proof}

\subsection{An alternative proof of the specialization theorem}

We now establish the equivalence between equality of Hilbert functions and membership in a marked family over a local coefficient ring.
We then use the weakly reverse lexicographic property and base change to obtain the specialization of reduced Gr\"obner bases.

Throughout this subsection and the next one, let $(B,\mathfrak m)$ be a local domain, and put $L:=\operatorname{Frac}(B)$ and $k:=B/\mathfrak m$.
We begin with an elementary local freeness criterion.

\begin{lemma}[Local freeness criterion]\label{lem:localfree}
Let $M$ be a finitely generated $B$-module.
Suppose that
\[
\dim_L(M\otimes_B L)
=
\dim_k(M/\mathfrak m M)
=
r.
\]
Then $M$ is a free $B$-module of rank $r$.
\end{lemma}

\begin{proof}
Choose a $k$-basis $\{\overline m_1,\ldots,\overline m_r\}$ of $M/\mathfrak mM$, and lift its elements to $m_1,\ldots,m_r\in M$.
By Nakayama's lemma, these elements generate $M$.
Hence there is a surjective $B$-linear map
\[
\pi:B^r\longrightarrow M \ ;\  e_i\longmapsto m_i.
\]
Here $e_1,\ldots,e_r$ are the standard basis vectors of $B^r$.
Put $Q:=\ker\pi$.
Since $L$ is a localization of $B$ by $S:=B\smallsetminus \{0\}$, say $L=S^{-1}B$, it is flat over $B$.
Hence the sequence $0\to Q\to B^r\xrightarrow{\pi}M\to0$ remains exact after tensoring with $L$ over $B$:
\[
0\longrightarrow Q\otimes_B L
\longrightarrow L^r
\xrightarrow{\pi_L}
M\otimes_B L
\longrightarrow0.
\]
By the assumption, the two vector spaces $L^r$ and $M\otimes_BL$ connected by $\pi_L$ have the same dimension $r$.
Therefore $\pi_L$ is an isomorphism, and hence $Q\otimes_B L=0$.

Since $Q\subset B^r$ and $B$ is an integral domain, $Q$ is torsion-free.
Consequently, the natural localization map $Q\to Q\otimes_B L \cong S^{-1}Q$ is injective.
It follows that $Q=0$, so $\pi$ is an isomorphism.
Thus we obtain $M\cong B^r$.
\end{proof}

For the remainder of this subsection and the next one, assume in addition that $B$ is Noetherian.
Write
\[
S_B:=B[x_1,\ldots,x_n], \qquad S_L:=L[x_1,\ldots,x_n], \qquad S_k:=k[x_1,\ldots,x_n].
\]
Let $I_B\subset S_B$ be a homogeneous ideal, and set $I_L:=I_BS_L$ and $I_k:=I_BS_k$.

\begin{lemma}\label{lem:graded-local-free}
For $d\geq0$, suppose that
\begin{equation}\label{eq:Hilbeq}
\HF_{S_L/I_L}(d) = \HF_{S_k/I_k}(d) =:r_d.
\end{equation}
Then the following hold:
\begin{enumerate}[label=\textup{(\roman*)}]
\item $(S_B/I_B)_d$ is a free $B$-module of rank $r_d$.
\item Let $J:=\ini_\prec(I_k)$.
The natural $B$-linear map
\[
\theta_d:
\langle\mathcal N(J)_d\rangle_B \longrightarrow (S_B/I_B)_d \ ; \  x^\beta\longmapsto x^\beta+I_B
\]
is an isomorphism.
In particular, the residue classes of $\mathcal N(J)_d$ form a $B$-basis of $(S_B/I_B)_d$.
\end{enumerate}
Hence, if \eqref{eq:Hilbeq} holds for every $d\geq0$, then $S_B=I_B\oplus_B\langle\mathcal N(J)\rangle_B$, that is, $I_B\in\mathbf{Mf}_J(B)$.
\end{lemma}

\begin{proof}
(i) Put $M_d:=(S_B/I_B)_d$.
Since $(S_B)_d$ is a free $B$-module of finite rank and $M_d=(S_B)_d/(I_B)_d$, the module $M_d$ is finitely generated over $B$.
Taking the homogeneous component of degree $d$ commutes with tensoring over $B$ with $L$ or $k$.
Thus $M_d\otimes_B L\cong(S_L/I_L)_d$ and $M_d\otimes_B k\cong(S_k/I_k)_d$ as $L$-linear spaces and $k$-linear spaces, respectively.
Since $M_d\otimes_B k\cong M_d/\mathfrak mM_d$ as $k$-linear spaces, the assumption \eqref{eq:Hilbeq} implies $\dim_L(M_d\otimes_B L) = \dim_k(M_d/\mathfrak mM_d) = r_d$.
Applying Lemma~\ref{lem:localfree}, we have $M_d\cong B^{r_d}$.

\medskip
\noindent
(ii) Since $J=\ini_\prec(I_k)$, Macaulay's basis theorem~\cite[Theorem~15.3]{Eisenbud} shows that the residue classes of $\mathcal N(J)_d$ form a $k$-basis of $(S_k/I_k)_d$.
Thus $\#\mathcal N(J)_d=\HF_{S_k/I_k}(d)=r_d$.
Together with (i), this shows that the source and target of $\theta_d$ are free $B$-modules of rank $r_d$.
Moreover, tensoring $\theta_d$ with $k$ over $B$ gives the map $\theta_d\otimes_B k:\langle\mathcal N(J)_d\rangle_k\longrightarrow(S_k/I_k)_d$, where $\theta_d\otimes_B k$ maps each $x^\beta$ to $x^\beta+I_k$.
Hence $\theta_d\otimes_B k$ is an isomorphism.

Recall that a homomorphism between finite free modules of the same rank over a local ring is an isomorphism if its reduction modulo the maximal ideal is an isomorphism.
Indeed, after choosing bases, its representing matrix $A$ is square, and invertibility modulo $\mathfrak m$ implies $\det(A)\notin\mathfrak m$.
Since $B$ is local, $\det(A)$ is a unit in $B$, so $A$ is invertible.
Applying this fact to $\theta_d$ and $\theta_d\otimes_B k$, we conclude that $\theta_d$ is an isomorphism.

\medskip

Finally, if \eqref{eq:Hilbeq} holds in every degree,
(ii) gives $(S_B)_d = (I_B)_d \oplus_B \langle\mathcal N(J)_d\rangle_B$ for any $d\geq 0$.
Taking the direct sum over all degrees yields $S_B=I_B\oplus_B\langle\mathcal N(J)\rangle_B$.
\end{proof}

We next clarify the relation between marked bases and reduced Gr\"obner bases, which is discussed in~\cite[Section~5]{LR16}.
The following lemma shows that a $J$-marked basis is the reduced Gr\"obner basis of the ideal it generates whenever $J$ is weakly reverse lexicographic.

\begin{lemma}[Marked bases and reduced Gr\"obner bases]\label{lem:wrl-marked-gb}
Let $\prec$ be the graded reverse lexicographic order with $x_1\succ\cdots\succ x_n$.
Let $J$ be a monomial ideal, let $\mathcal F_J=\{f_\alpha\mid x^\alpha\in\mathcal B_J\}$ be a $J$-marked basis over $C$, and put $I:=\langle\mathcal F_J\rangle$.
\begin{enumerate}[label=\textup{(\roman*)}]
\item If $\LM_\prec(f_\alpha)=x^\alpha$ for every $x^\alpha\in\mathcal B_J$, then $\mathcal F_J$ is the reduced Gr\"obner basis of $I$, and $\ini_\prec(I)=J$.
\item If $J$ is weakly reverse lexicographic with respect to $\prec$, then $\LM_\prec(f_\alpha)=x^\alpha$ for every $x^\alpha\in\mathcal B_J$.
Consequently, $\mathcal F_J$ is the reduced Gr\"obner basis of $I$.
\end{enumerate}
\end{lemma}

\begin{proof}
(i) Since the $f_\alpha$ are monic and their leading monomials are $x^\alpha$, ordinary polynomial reduction by $\mathcal F_J$ is available over $C$.
Each reduction decreases the monomials with respect to $\prec$, so the process terminates.

Here we show that every element of $I$ reduces to zero with respect to $\mathcal F_J$.
Let $f\in I$, and let $r$ be a remainder obtained by reduction with respect to $\mathcal F_J$.
Since $\LM_{\prec}(\mathcal{F}_J)=\mathcal{B}_J$, we have $r\in\langle\mathcal N(J)\rangle_C$.
On the other hand, $f-r\in I$, so $r\in I$.
As $\mathcal F_J$ is a $J$-marked basis, we have $I\cap\langle\mathcal N(J)\rangle_C=\{0\}$, which implies $r=0$.
Thus every element of $I$ reduces to zero with respect to $\mathcal F_J$, proving that $\mathcal F_J$ is a monic Gr\"obner basis and $\ini_\prec(I)=J$.

Moreover, the leading monomials of $\mathcal F_J$ are precisely the minimal monomial generators of $J$, and the tail of each $f_\alpha$ is supported on $\mathcal N(J)$.
Therefore $\mathcal F_J$ is reduced.

\medskip
\noindent
(ii) For each $x^\alpha\in\mathcal B_J$, write $f_\alpha$ as in \eqref{eq:marked-poly}.
If $x^\beta\in\mathcal N(J)_{|\alpha|}$ satisfied $x^\beta\succ x^\alpha$, the weakly reverse lexicographic property would imply $x^\beta\in J$, a contradiction.
Thus every tail monomial of $f_\alpha$ is smaller than $x^\alpha$, and hence $\LM_\prec(f_\alpha)=x^\alpha$.
The conclusion now follows from (i).
\end{proof}

We can now establish the specialization theorem over a local coefficient ring.

\begin{theorem}[Specialization via marked families]\label{thm:marked-specialization}
Let $\prec$ be the graded reverse lexicographic order with $x_1\succ\cdots\succ x_n$, and let $J:=\ini_\prec(I_k)$.
Then $\HS_{S_L/I_L}(z)=\HS_{S_k/I_k}(z)$ if and only if $S_B=I_B\oplus_B\langle\mathcal N(J)\rangle_B$, that is, $I_B\in\mathbf{Mf}_J(B)$.

If these equivalent conditions hold and $J$ is weakly reverse lexicographic, then $I_B$ has a unique $J$-marked basis $\mathcal F_{J,B}=\{f_{\alpha,B}\mid x^\alpha\in\mathcal B_J\}$, which is its reduced Gr\"obner basis with respect to $\prec$.
The base changes $\mathcal F_{J,L}$ and $\mathcal F_{J,k}$ obtained from $\mathcal F_{J,B}$ by applying the natural maps $B\to L$ and $B\to k$ coefficientwise, respectively, are the reduced Gr\"obner bases of $I_L$ and $I_k$, with the same leading monomials $\mathcal B_J$.
In particular, $\GB_\prec(I_L)\subset S_B$, and reducing its elements coefficientwise modulo $\mathfrak m$ gives $\GB_\prec(I_k)$.
Consequently, the initial ideals of $I_L$ and $I_k$ have the same minimal monomial generators $\mathcal B_J$, and both are weakly reverse lexicographic.
\end{theorem}

\begin{proof}
Suppose first that $\HS_{S_L/I_L}(z)=\HS_{S_k/I_k}(z)$.
Then Lemma~\ref{lem:graded-local-free} gives $S_B=I_B\oplus_B\langle\mathcal N(J)\rangle_B$.
Conversely, suppose that $S_B=I_B\oplus_B\langle\mathcal N(J)\rangle_B$.
Then $S_B/I_B\cong\langle\mathcal N(J)\rangle_B$ as graded $B$-modules.
Tensoring this isomorphism over $B$ with $L$ and $k$ gives $S_L/I_L\cong\langle\mathcal N(J)\rangle_L$ and $S_k/I_k\cong\langle\mathcal N(J)\rangle_k$, respectively.
Hence both fibers have the same Hilbert function $d\mapsto\#\mathcal N(J)_d$.

Assume that the equivalent conditions hold. 
If $J$ is weakly reverse lexicographic, then it is strongly stable, and hence Proposition~\ref{prop:marked-basis-I} gives the unique $J$-marked basis of $I_B$.
Since $J$ is weakly reverse lexicographic, Lemma~\ref{lem:wrl-marked-gb}(ii) shows that $\mathcal F_{J,B}$ is the reduced Gr\"obner basis of $I_B$.
By Corollary~\ref{cor:marked-BLk}, $\mathcal F_{J,L}$ and $\mathcal F_{J,k}$ are the unique $J$-marked bases of $I_L$ and $I_k$, respectively.
Applying Lemma~\ref{lem:wrl-marked-gb}(ii) again, we conclude that these are the reduced Gr\"obner bases of the two fibers, with the same leading monomials $\mathcal B_J$.
The assertion about coefficients and their specialization follows from the coefficientwise base-change description in Corollary~\ref{cor:marked-BLk}.
\end{proof}

We apply this result to the universal family introduced in Section~\ref{sec:prelim}.

\begin{proof}[Alternative proof of Corollary~\ref{cor:fixed-full-specialization}]
Let $\bm F\in V_\tau$ satisfy $\HS_{R/I_{\bm F}}(z)=H_\tau^{\mathrm{gen}}(z)$, and suppose that $\ini_\prec(I_{\bm F})$ is weakly reverse lexicographic.

Put $B:=\mathcal O_{\bm F}$ and consider $I_B:=\langle g_1,\ldots,g_m\rangle\subset B[x]$.
As established in Section~\ref{sec:prelim}, $B$ is a Noetherian local domain with fraction field $L_\tau$ and residue field $K$.
For the ideal $I_B$ above, we have $I_B L_\tau[x]=I_\tau^{\mathrm{gen}}$ and $\pi_{\bm F}(I_B)=I_{\bm F}$.

The assumed equality of Hilbert series gives $\HF_{L_\tau[x]/I_\tau^{\mathrm{gen}}}(d) = \HF_{R/I_{\bm F}}(d)$ for every $d \geq 0$.
Thus Theorem~\ref{thm:marked-specialization} applies with $J:=\ini_\prec(I_{\bm F})$.
It follows that $\GB_\prec(I_\tau^{\mathrm{gen}})\subset\mathcal O_{\bm F}[x]$ and $\pi_{\bm F}\bigl(\GB_\prec(I_\tau^{\mathrm{gen}})\bigr) =\GB_\prec(I_{\bm F})$, with the same leading monomials.
In particular, $\ini_\prec(I_\tau^{\mathrm{gen}})$ is weakly reverse lexicographic.
\end{proof}

\subsection{Specialization up to a prescribed degree by truncation}

The preceding specialization theorem assumes equality of Hilbert functions in every degree.
To obtain the corresponding result up to degree $D$, we truncate the ideals by adjoining $\langle x_1,\ldots,x_n\rangle^{D+1}$.
This preserves their homogeneous components through degree $D$ and makes the Hilbert functions of the truncated quotients vanish in higher degrees.
We can therefore apply Theorem~\ref{thm:marked-specialization} to the truncated family.

\begin{theorem}[Specialization up to a prescribed degree via marked families]\label{thm:marked-specialization-finite}
Let $\prec$ be the graded reverse lexicographic order with $x_1\succ\cdots\succ x_n$.
Let $D$ be a nonnegative integer, and assume that
\[
\HF_{S_L/I_L}(d)=\HF_{S_k/I_k}(d)
\qquad(0\leq d\leq D).
\]
Suppose that $J:=\ini_\prec(I_k)$ is weakly reverse lexicographic up to degree $D$.
Then, for every $0\leq d\leq D$, one has
$\GB_\prec(I_L)_d\subset S_B$, and
\[
\pi_k\bigl(\GB_\prec(I_L)_d\bigr)=\GB_\prec(I_k)_d,
\qquad
\LM_\prec\bigl(\GB_\prec(I_L)_d\bigr)
=\LM_\prec\bigl(\GB_\prec(I_k)_d\bigr),
\]
where $\pi_k:B\to k$ is the residue homomorphism, extended coefficientwise to polynomials over $B$.
In particular, $\ini_\prec(I_L)$ and $\ini_\prec(I_k)$ contain the same monomials of degree at most $D$.
\end{theorem}

\begin{proof}
For $C=B,L,k$, put
\[
\mathfrak n_C:=\langle x_1,\ldots,x_n \rangle\subset S_C,
\qquad
\widehat I_C:=I_C+\mathfrak n_C^{D+1}.
\]
The ideals $\widehat I_L$ and $\widehat I_k$ are the generic and special fibers of $\widehat I_B$, respectively.
Moreover,
\[
(S_C/\widehat I_C)_d\cong
\begin{cases}
(S_C/I_C)_d & (0\leq d\leq D),\\
0 & (d>D).
\end{cases}
\]
Hence $\HF_{S_L/\widehat I_L}(d) =\HF_{S_k/\widehat I_k}(d)$ for every $d\geq0$.

Since the ideals are homogeneous, we have
$\widehat J:=\ini_\prec(\widehat I_k)
=J+\mathfrak n_k^{D+1}$.
Every minimal monomial generator of $\widehat J$ has degree at most $D+1$.
Those of degree at most $D$ are minimal monomial generators of $J$ and satisfy the weakly reverse lexicographic condition by hypothesis.
For the remaining minimal generators, the condition holds because every monomial of degree $D+1$ belongs to $\widehat J$.
Thus $\widehat J$ is weakly reverse lexicographic, and Theorem~\ref{thm:marked-specialization} applies to $\widehat I_B$.
It gives the specialization of the reduced Gr\"obner bases of $\widehat I_L$ and $\widehat I_k$, with their leading monomials unchanged.

For $C=L,k$, the homogeneous components of $I_C$ and $\widehat I_C$ agree through degree $D$.
Therefore their initial ideals contain the same monomials of degree at most $D$, and their reduced Gr\"obner basis elements of degree at most $D$ coincide:
\[
\GB_\prec(\widehat I_C)_d=\GB_\prec(I_C)_d
\qquad(0\leq d\leq D).
\]
Indeed, in each such degree, the minimal leading monomials and the standard monomials are the same, so the corresponding reduced Gr\"obner basis elements are uniquely determined by the same homogeneous component of the ideal.
The assertions now follow from these equalities and Theorem~\ref{thm:marked-specialization} applied to $\widehat I_B$.
\end{proof}

\begin{proof}[Alternative proof of Theorem~\ref{thm:fixed-degreewise-specialization}]
The assertion follows directly from Theorem~\ref{thm:marked-specialization-finite} applied to $B$, $I_B$, and $D$ as in the statement.
\end{proof}

\subsection{Sequences with a prescribed projective zero}

Theorem~\ref{thm:marked-specialization} does not require the coefficients of the input polynomials to vary independently.
We apply it to the universal family on $Z_{\bm p}$ introduced in Section~\ref{sec:dimone}.
The coefficient space $Z_{\bm p}$ has a polynomial coordinate ring $A_\tau^{\bm p}$.
Localizing this ring at a specialization $\bm F\in Z_{\bm p}$ therefore places the restricted generic family in exactly the same setting as the preceding specialization theorem.
This also gives an alternative proof of the one-example criterion on $Z_{\bm p}$ without a change of coordinates.

For $\bm F\in Z_{\bm p}$, let $\pi_{\bm F}:A_\tau^{\bm p}\to K$ be the specialization defined in Section~\ref{sec:dimone}.
As in the proof of Theorem~\ref{thm:Zp-generic-HS}, put $\mathfrak n:=\ker\pi_{\bm F}$ and $B:=(A_\tau^{\bm p})_{\mathfrak n}$.
Then $\mathfrak n$ is a maximal ideal, and $B$ is a Noetherian local domain with fraction field $L_\tau^{\bm p}$ and residue field $K$.
We also denote by $\pi_{\bm F}$ the induced residue homomorphism $B\to K$ and its coefficientwise extension $B[x]\to K[x]$.

\begin{proof}[Alternative proof of the implication
\textup{(iii)} $\Rightarrow$ \textup{(i)} in
Theorem~\ref{thm:Zp-generic-HS}]
Let $\bm F\in Z_{\bm p}$ satisfy $\mathcal H_{\tau,\bm p}^{\mathrm{gen}}({\bm F})$ and $\mathcal W({\bm F})$, and put $I_B:=\langle\bm G_\tau^{\bm p}\rangle \subset B[x]$.
Its generic and special fibers are $\langle\bm G_\tau^{\bm p}\rangle\subset L_\tau^{\bm p}[x]$ and $I_{\bm F}\subset K[x]$, respectively.
By $\mathcal H_{\tau,\bm p}^{\mathrm{gen}}({\bm F})$,
the corresponding quotients have the same Hilbert function in every degree, and by $\mathcal W({\bm F})$, $J:=\ini_\prec(I_{\bm F})$ is weakly reverse lexicographic.
Thus Theorem~\ref{thm:marked-specialization} applies.
It shows that the reduced Gr\"obner basis of $\langle\bm G_\tau^{\bm p}\rangle$ has coefficients in $B$ and specializes to $\GB_\prec(I_{\bm F})$, with the same leading monomials.
Consequently, the initial ideal of the generic sequence is weakly reverse lexicographic.
\end{proof}

Together with Proposition~\ref{prop:H1-one-example},
this also recovers the implication \textup{(iii)}
$\Rightarrow$ \textup{(i)} in
Theorem~\ref{thm:dimone-criterion}.

\section{Magma code for the computational examples}\label{app:magma}

We give \textsc{Magma} code for searching for and verifying the examples in Proposition~\ref{prop:computed-types}.
The computations are performed over $\mathbb F_{31}$ in eight variables with graded reverse lexicographic order.
For $m\geq n$, the untruncated expression defining $H_\tau^{\mathrm{Fr}}(z)$ can be computed as an integer polynomial:
\[
\frac{\prod_{i=1}^m(1-z^{d_i})}{(1-z)^n}
=
\prod_{i=1}^n\left(\sum_{j=0}^{d_i-1}z^j\right)
\prod_{i=n+1}^m(1-z^{d_i}).
\]
The function \texttt{FrobergSeries} truncates this polynomial immediately before its first nonpositive coefficient.
The resulting polynomial is compared exactly with the Hilbert series returned by \texttt{HilbertSeries}.

The function \texttt{IsWRLFromGB} checks Definition~\ref{def:wrl} using the minimal monomial generators of the initial ideal.
In each degree containing such generators, it suffices to consider the monomials larger than the smallest generator of that degree.

\begin{lstlisting}[style=magmacode]
FrobergSeries := function(n, d)
    assert n ge 1 and #d ge n;
    assert &and[a ge 1 : a in d];

    P<t> := PolynomialRing(Integers());
    raw := &*[ &+[t^j : j in [0..d[i]-1]] : i in [1..n] ];
    for i in [n+1..#d] do
        raw *:= 1-t^d[i];
    end for;

    H := P!0;
    e := 0;
    while Coefficient(raw, e) gt 0 do
        H +:= Coefficient(raw, e)*t^e;
        e +:= 1;
    end while;
    return H;
end function;

IsWRLFromGB := function(G)
    LG := [LeadingMonomial(g) : g in G | g ne 0];
    if #LG eq 0 then
        return true;
    end if;
    R := Parent(LG[1]);

    for e in {Degree(u) : u in LG} do
        lowest := Min([u : u in LG | Degree(u) eq e]);
        LGsmall := [u : u in LG | Degree(u) le e];
        for t in MonomialsOfDegree(R, e) do
            if t gt lowest then
                contained := false;
                for s in LGsmall do
                    divides, _ := IsDivisibleBy(t, s);
                    if divides then
                        contained := true;
                        break;
                    end if;
                end for;
                if not contained then
                    return false;
                end if;
            end if;
        end for;
    end for;
    return true;
end function;

CheckExample := function(F, d)
    assert #F gt 0 and #F eq #d;
    R := Universe(F);
    for j in [1..#F] do
        if F[j] ne 0 then
            assert IsHomogeneous(F[j])
                   and Degree(F[j]) eq d[j];
        end if;
    end for;

    expected := FrobergSeries(Rank(R), d);
    I := ideal<R | F>;
    G := [g : g in GroebnerBasis(I) | g ne 0];
    H := HilbertSeries(I);
    N := Numerator(H);
    hasH := N eq (Parent(N)!expected)*Denominator(H);
    return hasH, IsWRLFromGB(G), H, G;
end function;

// Search parameters.
p := 31; n := 8; MaxTrials := 100; seed := 20261031;

K := GF(p);
R<[x]> := PolynomialRing(K, n, "grevlex");
Mons := [t : t in MonomialsOfDegree(R, 2)];
Sort(~Mons);
Reverse(~Mons);

v1, v2, v3 := GetVersion();
printf "Magma V%o.%o-%o\n", v1, v2, v3;
printf "p=%o; n=%o; base seed=%o\n", p, n, seed;
printf "Monomial exponents: %o\n",
       [Exponents(t) : t in Mons];

for m in [n..2*n] do
    SetSeed(seed+m);
    d := [2 : j in [1..m]];
    found := false;

    for trial in [1..MaxTrials] do
        C := [[Random(0, p-1) : t in Mons] : j in [1..m]];
        F := [&+[(K!C[j][i])*Mons[i] : i in [1..#Mons]] : j in [1..m]];
        hasH, hasW, H, G := CheckExample(F, d);
        if hasH and hasW then
            printf "m=%o; trial=%o; seed=%o\n",
                   m, trial, seed+m;
            printf "C_%o := %o;\n", m, C;
            printf "Hilbert series: %o\n", H;
            printf "Leading monomials: %o\n",
                   [LeadingMonomial(g) : g in G];
            found := true;
            break;
        end if;
    end for;

    if not found then
        printf "m=%o: no example found in %o trials.\n",
               m, MaxTrials;
    end if;
end for;
\end{lstlisting}

For each $m$, the search stops at the first example satisfying both conditions.
The rows of the printed integer matrix \texttt{C\_m} are the coefficient vectors of the input polynomials with respect to the ordered monomial sequence \texttt{Mons}.
These matrices, together with the monomial order and the characteristic, specify the examples independently of the random search.
Failure to find an example within the prescribed number of trials does not prove that no such example exists.

A recorded example can be verified without repeating the random search.
For example, to verify the case $m=8$, first define the three functions above and initialize \texttt{K}, \texttt{R}, and the ordered sequence \texttt{Mons} as in the search code.
Then execute the recorded assignment to \texttt{C\_8} and run:
\begin{lstlisting}[style=magmacode]
C := C_8;
m := #C;
F := [&+[(K!C[j][i])*Mons[i] : i in [1..#Mons]]
      : j in [1..m]];
hasH, hasW, H, G := CheckExample(F, [2 : j in [1..m]]);
assert hasH and hasW;
\end{lstlisting}
For the other values of $m$, use the corresponding coefficient matrix in the first line.

The code verifies both conditions over $\mathbb F_{31}$. 
Corollary~\ref{cor:arithmetic-lifting} then shows that both conjectures hold for these types over every field of characteristic zero.

\end{document}